\documentclass[leqno]{article}
\usepackage{amsmath,amsthm,amsfonts,amssymb}
\usepackage{mathrsfs}

\allowdisplaybreaks

\newfam\msbfam
\font\tenmsb=msbm10    \textfont\msbfam=\tenmsb \font\sevenmsb=msbm7
\scriptfont\msbfam=\sevenmsb \font\fivemsb=msbm5
\scriptscriptfont\msbfam=\fivemsb

\newfam\bigfam
\font\tenbig=msbm10 scaled \magstep2   \textfont\bigfam=\tenbig
\font\sevenbig=msbm7 scaled \magstep2 \scriptfont\bigfam=\sevenbig
\font\fivebig=msbm5 scaled \magstep2
\scriptscriptfont\bigfam=\fivebig

\begin{document}

\title{{\bf Anisotropic grand Herz type spaces with variable exponents and their applications}
 \footnotetext{\noindent *Corresponding author.}
}
\author{Hongbin WANG{$^{*}$} and Zongguang LIU}
\date{}
\maketitle

\begin{minipage}{13.5cm}
{\bf Abstract}\quad  \small{In this paper, we
introduce some anisotropic grand Herz type spaces with variable exponents, including anisotropic grand Herz spaces, anisotropic grand Herz-Morrey spaces and anisotropic grand Herz-type Hardy spaces with variable exponents. We obtain some properties and characterizations of these spaces in terms of some decompositions. Using their decompositions, we obtain some boundedness on the anisotropic grand Herz type spaces with variable exponents for some singular integral operators.}

\medskip


\medskip

{\bf MR(2020)Subject Classification}\quad 42B35, 42B20.

\end{minipage}

\section{ Introduction\label{s1}}$\indent$

The main purpose of this article is to introduce several types of anisotropic grand Herz type spaces with variable exponents and give the boundedness for a class of operators on these spaces.

The classical Herz type spaces mainly include Herz spaces, Herz-Morrey spaces and Herz-type Hardy spaces. The Herz spaces are a class of function spaces introduced by Herz [16] in the study of absolutely convergent Fourier transforms in 1968. The Herz-Morrey spaces can be traced to the work of Lu and Xu [27]. The complete theory of Herz-type Hardy spaces was established by Lu and Yang [28,29]. The classical Herz type spaces and their related function spaces have been playing an important role in harmonic analysis including the boundedness of many operators and partial differential equations, see [14,21,24,33] for example.

In recent years, as a generalization of the classical $L^p$ spaces, the Lebesgue spaces with variable exponent were introduced by Orlicz [31] in 1931. After the paper of Kov\'{a}\v{c}ik and R\'{a}kosn\'{i}k [22] in 1991, the theory of function spaces with variable exponents has attracted a lot of attentions in recent years(see [7,8] and the references therein). Owing to their wide applications in electrorheological fluids [34], image processing [6] and partial differential equations with non-standard growth [15], such theory has made great progresses. Particularly, the related theory of Herz type spaces with variable exponent has begun to be studied by many authors. For example, Izuki [18,19] considered the Herz spaces with variable exponent $\dot{K}^{\alpha,p}_{q(\cdot)}(\mathbb{R}^{n})$ ($K^{\alpha,p}_{q(\cdot)}(\mathbb{R}^{n})$) and the Herz-Morrey spaces with variable exponent $M\dot{K}^{\alpha,\lambda}_{q(\cdot),p}(\mathbb{R}^{n})$ in 2010. Subsequently, a certain Herz-type Hardy spaces with variable exponent $H\dot{K}^{\alpha,p}_{q(\cdot)}(\mathbb{R}^{n})$ ($HK^{\alpha,p}_{q(\cdot)}(\mathbb{R}^{n})$) were introduced by Wang and Liu [37] in 2012. Almeida and Drihem [1] introduced the Herz spaces with two variable exponents $\dot{K}^{\alpha(\cdot),p}_{q(\cdot)}(\mathbb{R}^{n})$ ($K^{\alpha(\cdot),p}_{q(\cdot)}(\mathbb{R}^{n})$). Furthermore, the Herz-type Hardy spaces with three variable exponents $H\dot{K}^{\alpha(\cdot),p(\cdot)}_{q(\cdot)}(\mathbb{R}^{n})$ ($HK^{\alpha(\cdot),p(\cdot)}_{q(\cdot)}(\mathbb{R}^{n})$) were introduced by Drihem and Seghiri [10] in 2016. Meanwhile, the boundedness of some operators on above spaces have also been proved.

On the other hand, extending classical function spaces arising in harmonic analysis of Euclidean spaces to other domains and non-isotropic settings is an important topic. In 2003, Bownik [2] introduced the anisotropic Hardy spaces associated with very general discrete groups of dilations. The above spaces include the classical isotropic Hardy space theory of Fefferman and Stein [11] and parabolic Hardy space theory of Calder\'{o}n and Torchinsky [3,4]. In 2006, Lan [23] defined the anisotropic Herz spaces $\dot{K}^{\alpha,p}_{q}(A;\mathbb{R}^{n})$ ($K^{\alpha,p}_{q}(A;\mathbb{R}^{n})$) associated with the dilation $A$. A class of anisotropic Herz-type Hardy spaces associated with a non-isotropic dilation $H\dot{K}^{\alpha,p}_{q}(A;\mathbb{R}^{n})$ ($HK^{\alpha,p}_{q}(A;\mathbb{R}^{n})$) by Ding, Lan and Lu [9] in 2008. Liu, Yang and Yuan [25] introduced the anisotropic Hardy-Lorentz spaces $H^{p,q}_{A}(\mathbb{R}^{n})$ in 2016.

Combining the theory of variable exponent function spaces and anisotropic function spaces, Wang [36] defined the anisotropic Herz spaces with two variable exponents $\dot{K}^{\alpha(\cdot),p}_{q(\cdot)}(A;\mathbb{R}^{n})$ ($K^{\alpha(\cdot),p}_{q(\cdot)}(A;\mathbb{R}^{n})$) in 2015. Wang and Wu [38] introduced the anisotropic Herz-Morrey spaces with two variable exponents $M\dot{K}^{\alpha(\cdot),\lambda}_{q(\cdot),p}(A;\mathbb{R}^{n})$ ($MK^{\alpha(\cdot),\lambda}_{q(\cdot),p}(A;\mathbb{R}^{n})$) in 2016. Liu, Yang and Yuan [26] introduced the anisotropic variable Hardy-Lorentz space $H^{p(\cdot),q}_{A}(\mathbb{R}^{n})$ in 2017. Zhao and Zhou [39] gave the anisotropic Herz-type Hardy spaces with variable exponent $H\dot{K}^{\alpha,p}_{q(\cdot)}(A;\mathbb{R}^{n})$ ($HK^{\alpha,p}_{q(\cdot)}(A;\mathbb{R}^{n})$) in 2018.

The grand Lebesgue spaces were introduced by Iwaniec and Sbordone [17] to study the integrability of the Jacobian in 1992. The properties of
the grand Lebesgue spaces had been deeply investigated in [5,12], and the grand spaces have been proved to be useful in application to partial differential equations in [13,20]. The grand Lebesgue sequence spaces were introduced in [32], where various operators in harmonic analysis were studied on these spaces. Furthermore, Nafis, Rafeiro and Zaighum [30] introduced the grand variable Herz spaces $\dot{K}^{\alpha,p),\theta}_{q(\cdot)}(\mathbb{R}^{n})$ ($K^{\alpha,p),\theta}_{q(\cdot)}(\mathbb{R}^{n})$) in 2020. The grand Herz-Hardy spaces with variable exponent $H\dot{K}^{\alpha,p),\theta}_{q(\cdot)}(\mathbb{R}^{n})$ ($HK^{\alpha,p),\theta}_{q(\cdot)}(\mathbb{R}^{n})$) were introduced by Shabbir and Zaighum [35] in 2024.

Inspired by the above references, we introduce some anisotropic grand Herz type spaces with two variable exponents, including the anisotropic grand Herz spaces, the anisotropic grand Herz-Morrey spaces and the anisotropic grand Herz-type Hardy spaces with two variable exponents, which is a generalization of the anisotropic Herz type spaces and the Herz type spaces
with variable exponents. We obtain some properties and characterizations of these spaces in terms of some decompositions. Using their decompositions, we obtain some boundedness on the anisotropic grand Herz type spaces with two variable exponents for some operators.

To be precise, we explain the outline of this paper. We first briefly recall some standard notations and lemmas in Section 2. In Section 3, we will define the anisotropic Herz spaces with
two variable exponents $\dot{K}^{\alpha(\cdot),p}_{q(\cdot)}(A;\mathbb{R}^{n})$ and
$K^{\alpha(\cdot),p}_{q(\cdot)}(A;\mathbb{R}^{n})$, and give the
boundedness of some sublinear operators. Subsequently, the anisotropic Herz-Morrey spaces with
two variable exponents $M\dot{K}^{\alpha(\cdot),\lambda,\theta}_{q(\cdot),p)}(A;\mathbb{R}^{n})$ and
$MK^{\alpha(\cdot),\lambda,\theta}_{q(\cdot),p)}(A;\mathbb{R}^{n})$ will be defined and the
boundedness of some sublinear operators will be proved in Section 4. In Section 5, we will define the anisotropic Herz-type Hardy spaces with
two variable exponents $H\dot{K}^{\alpha(\cdot),p)}_{q(\cdot),\theta}(A;\mathbb{R}^{n})$ and
$HK^{\alpha(\cdot),p)}_{q(\cdot),\theta}(A;\mathbb{R}^{n})$. Moreover, the atomic decomposition theorems and the
boundedness of some linear operators on $H\dot{K}^{\alpha(\cdot),p)}_{q(\cdot),\theta}(A;\mathbb{R}^{n})$ and
$HK^{\alpha(\cdot),p)}_{q(\cdot),\theta}(A;\mathbb{R}^{n})$ will also be obtained.

In addition, we denote the Lebesgue measure and the characteristic
function of a measurable set $A\subset \mathbb{R}^{n}$ by $|A|$ and
$\chi_A$ respectively. $C$ always means a positive constant independent
of the main parameters and may change from one occurrence to another. The notation $f\lesssim g$ means that there exists a constant $C>0$ such that $f\leq Cg$. If $f\lesssim g$ and $g\lesssim f$, then $f\approx g$.

\section{ Preliminaries\label{s2}}$\indent$

We first introduce some basic definitions and properties of non-isotropic spaces associated with general expansive dilations. A $n\times n$ real matrix $A$ is called an expansive matrix, sometimes called a dilation, if all eigenvalues $\lambda$ of $A$ satisfy $|\lambda|>1$. We suppose $\lambda_1,..., \lambda_n$ are eigenvalues of $A$ (taken according to the multiplicity) so that $1<|\lambda_1|\leq...\leq|\lambda_n|$, and let $\lambda_-$, $\lambda_+$ be any numbers satisfying $1<\lambda_-<|\lambda_1|\leq |\lambda_n|<\lambda_+$. A set $\Delta\subset\mathbb{R}^n$ is said to be an ellipsoid if $\Delta=\{x\in \mathbb{R}^n: |Px|<1\}$, for some nondegenerate $n\times n$ matrix $P$, where $|\cdot|$ denotes the Euclidean norm in $\mathbb{R}^n$. For a dilation $A$, there exists an ellipsoid $\Delta$ and $r>1$ such that $\Delta\subset r\Delta\subset A\Delta$, where $|\Delta|$, the Lebesgue measure of $\Delta$, equals 1. Let $B_k=A^k\Delta$ for $k\in \mathbb{Z}$, then we have $B_k\subset rB_k\subset B_{k+1}$, and $|B_k|=b^k$, where $b=|\mathrm{det} A|=\prod_{j=1}^n|\lambda_j|$. Let $w$ be the smallest integer so that $2B_0\subset A^wB_0=B_w$. A quasi-norm associated with an expansive matrix $A$ is a measurable mapping $\rho_A: \mathbb{R}^n\rightarrow [0, \infty)$ satisfying
$$\rho_A(x)>0\hspace{3cm}\mathrm{for}\,\,x\neq 0,$$
$$\rho_A(Ax)=|\mathrm{det}A|\rho(x)\hspace{1.5cm}\mathrm{for}\,\,x\in \mathbb{R}^n,$$
$$\rho_A(x+y)\leq C(\rho_A(x)+\rho_A(y))\hspace{0.35cm}\mathrm{for}\,\,x,y\in \mathbb{R}^n,$$
where $C\geq 1$ is a constant. One can show that all quasi-norms associated to a fixed dilation $A$ are equivalent, see [2, Lemma 2.4]. Define the step homogeneous quasi-norm $\rho$ on $\mathbb{R}^n$ induced by dilation $A$ as
$$\displaystyle \rho(x)=\left\{\begin{array}{ll}
\displaystyle b^j\quad \mathrm{if}\,x\in B_{j+1}\setminus B_j,\\
\displaystyle 0,\quad \mathrm{if}\,x=0.
\end{array}\right.$$
For any $x,y\in \mathbb{R}^n$, we have
$$\rho(x+y)\leq b^w(\rho(x)+\rho(y)).\eqno(2.1)$$

Next we will recall some notations in variable function spaces. Given an open set $\Omega\subset \mathbb{R}^{n}$, and a measurable
function $p(\cdot):\Omega\longrightarrow(0,\infty),$
$L^{p(\cdot)}(\Omega)$ denotes the set of measurable functions $f$
on $\Omega$ such that for some $\lambda>0,$
$$\int_\Omega\left(\frac{|f(x)|}{\lambda}\right)^{p(x)}dx < \infty.$$
This set becomes a Banach function space when equipped with the
Luxemburg-Nakano norm
$$\|f\|_{L^{p(\cdot)}(\Omega)}=\inf\left\{\lambda>0:\int_\Omega
\left(\frac{|f(x)|}{\lambda}\right)^{p(x)}dx \leq 1\right\}.$$ These
spaces are referred to as variable $L^{p}$ spaces, since they generalized the standard
$L^{p}$ spaces.

The space $L_{\rm
loc}^{p(\cdot)}(\Omega)$ is defined by $L_{\rm
loc}^{p(\cdot)}(\Omega):=\{f: f\in L^{p(\cdot)}(E) \,\,\mathrm{for \,\,all \,\,compact \,\,subsets} \,\,E\subset \Omega\}.$
Define
$\mathcal{P}^0(\Omega)$ to be set of
$p(\cdot):\Omega\longrightarrow(0,\infty)$ such that
$$p^{-}=\mathrm{ess} \inf\{p(x):x\in \Omega\}>0,\quad  p^{+}=\mathrm{ess} \sup\{p(x):x\in \Omega\}<\infty.$$
Define $\mathcal{P}(\Omega)$ to be set of
$p(\cdot):\Omega\longrightarrow[1,\infty)$ such that
$$p^{-}=\mathrm{ess} \inf\{p(x):x\in\Omega\}>1,\quad  p^{+}=\mathrm{ess} \sup\{p(x):x\in\Omega\}<\infty.$$
Denote $p'(x)=p(x)/(p(x)-1).$

Let $f$ be a locally integrable function. The Hardy-Littlewood
maximal operator is defined by
$$Mf(x)=\sup_{B\ni x}\frac{1}{|B|}\int_{B\cap \Omega}|f(y)|dy,$$
where the supremum is taken over all balls $B$ containing $x$.
Let $\mathcal{B}(\Omega)$ be the set of $p(\cdot)\in$
$\mathcal{P}(\Omega)$
 such that the Hardy-Littlewood maximal operator
$M$ is bounded on $L^{p(\cdot)}(\Omega)$.

In variable $L^{p}$ spaces there are some important lemmas as
follows.

\noindent{\bf Lemma 2.1}$^{[22]}$\quad Let $p(\cdot)\in
\mathcal{P}(\mathbb{R}^{n})$. If $f\in L^{p(\cdot)}(\mathbb{R}^{n})$
and $g\in L^{p'(\cdot)}(\mathbb{R}^{n})$, then $fg$ is integrable on
$\mathbb{R}^{n}$ and
$$\int_{\mathbb{R}^{n}}|f(x)g(x)|dx \leq r_{p}\|f\|_{L^{p(\cdot)}}\|g\|_{L^{p'(\cdot)}},$$
where $$r_{p}=1+1/p^--1/p^+.$$

This inequality is named the generalized H\"{o}lder inequality with
respect to the variable $L^{p}$ spaces.

\noindent{\bf Lemma 2.2}$^{[18]}$\quad Suppose $p(\cdot)\in
\mathcal{B}(\mathbb{R}^{n})$. Then there exists a constant $C>0$
such that for all balls $B$ in $\mathbb{R}^{n}$,
$$\frac{1}{|B|}\|\chi_B\|_{L^{p(\cdot)}(\mathbb{R}^{n})}\|\chi_B\|_{L^{p'(\cdot)}(\mathbb{R}^{n})}\leq
C.$$

\noindent{\bf Lemma 2.3}$^{[18]}$\quad Let $p(\cdot)\in
\mathcal{B}(\mathbb{R}^{n})$. Then for all balls $B$ in $\mathbb{R}^{n}$ and all
measurable subsets $S\subset B$,
$$\frac{\|\chi_B\|_{L^{p(\cdot)}(\mathbb{R}^{n})}}{\|\chi_S\|_{L^{p(\cdot)}(\mathbb{R}^{n})}}\lesssim\frac{|B|}{|S|},\,\,
\frac{\|\chi_S\|_{L^{p(\cdot)}(\mathbb{R}^{n})}}{\|\chi_B\|_{L^{p(\cdot)}(\mathbb{R}^{n})}}\lesssim\left(\frac{|S|}{|B|}\right)^{\delta_1} \mathrm{and}\,\, \frac{\|\chi_S\|_{L^{p'(\cdot)}(\mathbb{R}^{n})}}{\|\chi_B\|_{L^{p'(\cdot)}(\mathbb{R}^{n})}}\lesssim\left(\frac{|S|}{|B|}\right)^{\delta_2},$$
where $0<\delta_1, \delta_2<1$ are constants.

\noindent {\bf Remark 2.1}\quad Throughout this paper $\delta_2$ is the same as in Lemma 2.3.

\noindent{\bf Lemma 2.4}$^{[7]}$\quad Let $p(\cdot),q(\cdot),r(\cdot)\in
\mathcal{P}(\mathbb{R}^{n})$ such that $1/p(x)=1/q(x)+1/r(x)$. If $f\in L^{q(\cdot)}(\mathbb{R}^{n})$
and $g\in L^{r(\cdot)}(\mathbb{R}^{n})$, then $fg\in L^{p(\cdot)}(\mathbb{R}^{n})$ and
$$\|fg\|_{L^{p(\cdot)}(\mathbb{R}^{n})} \leq C\|f\|_{L^{q(\cdot)}(\mathbb{R}^{n})}\|g\|_{L^{r(\cdot)}(\mathbb{R}^{n})},$$
where $C$ is a constant independent of the functions $f$ and $g$.

We can obtain the following definition in [1].

\noindent {\bf Definition 2.1}$^{[1]}$\quad Let a function $g(\cdot): \mathbb{R}^{n}\rightarrow\mathbb{R}$.

\noindent (1) $g(\cdot)$ is locally log-H\"{o}lder continuous, if there exists a constant $C>0$ such that
$$|g(x)-g(y)|\leq \frac{C}{\log(e+1/|x-y|)}$$
for all $x,y\in \mathbb{R}^{n}$ and $|x-y|<1/2$.

\noindent (2) $g(\cdot)$ is locally log-H\"{o}lder continuous at the origin (or has a log decay at the origin), if there exists a constant $C>0$ such that
$$|g(x)-g(0)|\leq \frac{C}{\log(e+1/|x|)}$$
for all $x\in \mathbb{R}^{n}$.

\noindent (3) $g(\cdot)$ is locally log-H\"{o}lder continuous at infinity (or has a log decay at infinity), if there exist some $g_\infty\in\mathbb{R}^{n}$ and $C>0$ such that
$$|g(x)-g_\infty|\leq \frac{C}{\log(e+|x|)}$$
for all $x\in \mathbb{R}^{n}$.

By $\mathcal{P}_0(\mathbb{R}^{n})$ and $\mathcal{P}_\infty(\mathbb{R}^{n})$ we denote the class of all exponents $p\in\mathcal{P}(\mathbb{R}^{n})$ which are locally log-H\"{o}lder continuous at the origin and at infinity, respectively.

We introduce the grand Lebesgue sequence spaces at the end of this section. Denote
$\mathbb{Z_+}$ and $\mathbb{N}$ as the sets of all positive and non-negative integers. The letter $\mathbb{X}$ stands for one of the sets $\mathbb{Z}^n$, $\mathbb{Z}$, $\mathbb{Z}_+$ and $\mathbb{N}$.

\noindent {\bf Definition 2.2}$^{[32]}$\quad Let $1\leq p<\infty$ and $\theta>0$. The grand Lebesgue sequence spaces
$l^{p),\theta}(\mathbb{X})$ is defined by the norm
$$\|X\|_{l^{p),\theta}(\mathbb{X})}:=\sup_{\varepsilon>0}\left(\varepsilon^\theta\sum_{k\in\mathbb{X}}|x_k|^{p(1+\varepsilon)}\right)^{\frac{1}{p(1+\varepsilon)}}=\sup_{\varepsilon>0}\varepsilon^{\frac{\theta}{p(1+\varepsilon)}}\|x\|^{l^{p(1+\varepsilon)}},$$
where $X=\{x_k\}_{k\in\mathbb{X}}$.

Note that the following nesting properties hold:
$$l^{p(1-\varepsilon)}\hookrightarrow l^p\hookrightarrow l^{p),\theta_1}\hookrightarrow l^{p),\theta_2}\hookrightarrow l^{p(1+\delta)},$$
for $0<\varepsilon<\frac{1}{p}$, $\delta>0$ and $0<\theta_1\leq \theta_2$.

\section{Anisotropic grand Herz spaces with variable
exponents\label{s3}}$\indent$

In this section, we first introduce the definition of anisotropic grand Herz spaces with two variable
exponents. Let $C_k=B_k\setminus B_{k-1}$ for $k\in \mathbb{Z}$, $\chi_k=\chi_{C_k}$ for $k\in \mathbb{Z}$,
$\tilde{\chi}_k=\chi_k$ if $k\in \mathbb{Z_+}$ and
$\tilde{\chi}_0=\chi_{B_0}$, where $\chi_{C_k}$ is the
characteristic function of $C_k$.

\noindent {\bf Definition 3.1}\quad Let $\alpha(\cdot):
\mathbb{R}^{n}\rightarrow\mathbb{R}$ with $\alpha(\cdot)\in L^\infty(\mathbb{R}^{n})$, $1\leq p<\infty$, $q(\cdot)\in
\mathcal{P}(\mathbb{R}^{n})$ and $\theta>0$. The homogeneous anisotropic grand Herz space
$\dot{K}^{\alpha(\cdot),p)}_{q(\cdot),\theta}(A;\mathbb{R}^{n})$ associated with the dilation $A$ is defined by
$$\dot{K}^{\alpha(\cdot),p)}_{q(\cdot),\theta}(A;\mathbb{R}^{n})=\{f\in L_{\rm
loc}^{q(\cdot)}(\mathbb{R}^{n}\setminus \{0\}):
\|f\|_{\dot{K}^{\alpha(\cdot),p)}_{q(\cdot),\theta}(A;\mathbb{R}^{n})}<\infty\},$$
where
$$\|f\|_{\dot{K}^{\alpha(\cdot),p)}_{q(\cdot),\theta}(A;\mathbb{R}^{n})}=\sup_{\varepsilon>0}\left(\varepsilon^\theta\sum_{k=-\infty}^\infty\|b^{k\alpha(\cdot)
}f\chi_k\|^{p(1+\varepsilon)}_{L^{q(\cdot)}(\mathbb{R}^{n})}\right)^{\frac{1}{p(1+\varepsilon)}}.$$ The
non-homogeneous anisotropic grand Herz space $K^{\alpha(\cdot),p)}_{q(\cdot),\theta}(A;\mathbb{R}^{n})$ associated with the dilation $A$
is defined by $$K^{\alpha(\cdot),p)}_{q(\cdot),\theta}(A;\mathbb{R}^{n})=\{f\in
L_{\rm loc}^{q(\cdot)}(\mathbb{R}^{n}):
\|f\|_{K^{\alpha(\cdot),p)}_{q(\cdot),\theta}(A;\mathbb{R}^{n})}<\infty\},$$ where
$$\|f\|_{K^{\alpha(\cdot),p)}_{q(\cdot),\theta}(A;\mathbb{R}^{n})}=\sup_{\varepsilon>0}\left(\varepsilon^\theta\sum_{k=0}^\infty\|b^{k\alpha(\cdot)}f\tilde{\chi}_k\|^{p(1+\varepsilon)}_{L^{q(\cdot)}(\mathbb{R}^{n})}\right)^{\frac{1}{p(1+\varepsilon)}}.$$

Similar to the method in [10, Proposition 1], we can obtain the following lemma.

\noindent{\bf Lemma 3.1}\quad Let $1\leq p<\infty$, $q(\cdot)\in
\mathcal{P}^0(\mathbb{R}^{n})$ and $\theta>0$. If $\alpha(\cdot)\in L^\infty(\mathbb{R}^{n})\cap\mathcal{P}_\infty(\mathbb{R}^{n})$, then
$$K^{\alpha(\cdot),p)}_{q(\cdot),\theta}(A;\mathbb{R}^{n})=K^{\alpha_\infty,p)}_{q(\cdot),\theta}(A;\mathbb{R}^{n}).$$
If $\alpha(\cdot)\in L^\infty(\mathbb{R}^{n})\cap\mathcal{P}_0(\mathbb{R}^{n})$, then
$$\|f\|_{\dot{K}^{\alpha(\cdot),p)}_{q(\cdot),\theta}(A;\mathbb{R}^{n})}\approx\sup_{\varepsilon>0}\varepsilon^{\frac{\theta}{p(1+\varepsilon)}}\left(\sum_{k=-\infty}^{-1}\|b^{k\alpha(0)
}f\chi_k\|^{p(1+\varepsilon)}_{L^{q(\cdot)}(\mathbb{R}^{n})}+\sum_{k=0}^\infty\|b^{k\alpha_\infty
}f\chi_k\|^{p(1+\varepsilon)}_{L^{q(\cdot)}(\mathbb{R}^{n})}\right)^{\frac{1}{p(1+\varepsilon)}}.$$

Next we will consider the decomposition of
$\dot{K}^{\alpha(\cdot),p)}_{q(\cdot),\theta}(A;\mathbb{R}^{n})$. We begin
with the notation of central block.

\noindent{\bf Definition 3.2}$^{[36]}$\quad Let $q(\cdot)\in
\mathcal{P}(\mathbb{R}^{n}), \alpha(\cdot)\in L^\infty(\mathbb{R}^{n})\cap\mathcal{P}_0(\mathbb{R}^{n})\cap\mathcal{P}_\infty(\mathbb{R}^{n})$ and $0<\alpha_l<\infty$. Denote $\alpha_l=\alpha(0), l<0; \alpha_l=\alpha_\infty, l\geq 0$.

(i) A measurable function $a(x)$ is said to be a central
$(\alpha(\cdot), q(\cdot))$-block if

\hspace{3mm}(1) supp\,$a\subset B_l$.

\hspace{3mm}(2) $\|a\|_{L^{q(\cdot)}(\mathbb{R}^{n})}\leq
b^{-l\alpha_l}$.

(ii) A measurable function $a(x)$ is said to be a central
$(\alpha(\cdot), q(\cdot))$-block of restricted type if

\hspace{3mm}(1) supp\,$a\subset B_l$ for some $l\geq 0$.

\hspace{3mm}(2) $\|a\|_{L^{q(\cdot)}(\mathbb{R}^{n})}\leq
b^{-l\alpha_\infty}$.

The following decomposition theorem shows that the central blocks
are the ``building block" of the anisotropic grand Herz spaces with exponents.

\noindent{\bf Theorem 3.1}\quad Let $1\leq p<\infty, q(\cdot)\in \mathcal{P}(\mathbb{R}^{n}), \alpha(\cdot)\in L^\infty(\mathbb{R}^{n})\cap\mathcal{P}_0(\mathbb{R}^{n})\cap\mathcal{P}_\infty(\mathbb{R}^{n})$
and $0<\alpha(0), \alpha_\infty<\infty$. The following two
statements are equivalent:

(i) $f\in \dot{K}^{\alpha(\cdot),p)}_{q(\cdot),\theta}(A;\mathbb{R}^{n})$.

(ii) $f$ can be represented by
$$f(x)=\sum_{k=-\infty}^\infty\lambda_kb_k(x),\eqno(3.1)$$ where each
$b_k$ is a central $(\alpha(\cdot), q(\cdot))$-block with support
contained in $B_k$ and $$\displaystyle\sup_{\varepsilon>0}\left(\varepsilon^\theta\sum_{k=-\infty}^\infty|\lambda_k|^{p(1+\varepsilon)}\right)^{\frac{1}{p(1+\varepsilon)}}<\infty.$$

\noindent{\bf Proof}\quad We first prove (i) implies (ii). For every
$f\in \dot{K}^{\alpha(\cdot),p)}_{q(\cdot),\theta}(A;\mathbb{R}^{n})$, write

$$\begin{array}{rl}
\displaystyle f(x)&\displaystyle=\sum_{k=-\infty}^\infty f(x)\chi_k(x)\\
&\displaystyle=\sum_{k=-\infty}^\infty\|b^{k\alpha(\cdot)}f\chi_k\|_{L^{q(\cdot)}(\mathbb{R}^{n})}\frac{f(x)\chi_k(x)}{\|b^{k\alpha(\cdot)}f\chi_k\|_{L^{q(\cdot)}(\mathbb{R}^{n})}}\\
&\displaystyle=\sum_{k=-\infty}^\infty\lambda_kb_k(x),

\end{array}$$
where
$\lambda_k=\|b^{k\alpha(\cdot)}f\chi_k\|_{L^{q(\cdot)}(\mathbb{R}^{n})}$
and
$b_k(x)=\frac{f(x)\chi_k(x)}{\|b^{k\alpha(\cdot)}f\chi_k\|_{L^{q(\cdot)}(\mathbb{R}^{n})}}$.

It is obvious that supp\,$b_k\subset B_k$ and
$\|b_k\|_{L^{q(\cdot)}(\mathbb{R}^{n})}=b^{-k\alpha_k}$. Thus,
each $b_k$ is a central $(\alpha(\cdot), q(\cdot))$-block with the
support $B_k$ and
$$\begin{array}{rl}
\displaystyle\sup_{\varepsilon>0}\left(\varepsilon^\theta\sum_{k=-\infty}^\infty|\lambda_k|^{p(1+\varepsilon)}\right)^{\frac{1}{p(1+\varepsilon)}}
&\displaystyle=\sup_{\varepsilon>0}\varepsilon^{\frac{\theta}{p(1+\varepsilon)}}\left(\sum_{k=-\infty}^\infty\|b^{k\alpha(\cdot)}f\chi_k\|^{p(1+\varepsilon)}_{L^{q(\cdot)}(\mathbb{R}^{n})}\right)^{\frac{1}{p(1+\varepsilon)}}\\
&\displaystyle=\sup_{\varepsilon>0}\varepsilon^{\frac{\theta}{p(1+\varepsilon)}}\|f\|_{\dot{K}^{\alpha(\cdot), p(1+\varepsilon)}_{q(\cdot)}(A;\mathbb{R}^{n})}\\
&\displaystyle=\|f\|_{\dot{K}^{\alpha(\cdot),p)}_{q(\cdot),\theta}(A;\mathbb{R}^{n})}<\infty.

\end{array}$$

Now we prove (ii) implies (i). Let $\displaystyle
f(x)=\sum_{k=-\infty}^\infty\lambda_kb_k(x)$ be a decomposition of $f$
which satisfies the hypothesis (ii) of Theorem 3.1. For each $j\in
\mathbb{Z}$, by the Minkowski inequality, we have
$$\|f\chi_j\|_{L^{q(\cdot)}(\mathbb{R}^{n})}\leq\sum_{k=j}^\infty|\lambda_k|\|b_k\|_{L^{q(\cdot)}(\mathbb{R}^{n})}.\eqno(3.2)$$

So we have
$$\begin{array}{rl}
\displaystyle\|f\|_{\dot{K}^{\alpha(\cdot),p)}_{q(\cdot),\theta}(A;\mathbb{R}^{n})}&\displaystyle=\sup_{\varepsilon>0}\left(\varepsilon^{\theta}\sum_{j=-\infty}^\infty\|b^{j\alpha(\cdot)}f\chi_j\|^{p(1+\varepsilon)}_{L^{q(\cdot)}(\mathbb{R}^{n})}\right)^{\frac{1}{p(1+\varepsilon)}}\\
&\displaystyle=\sup_{\varepsilon>0}\varepsilon^{\frac{\theta}{p(1+\varepsilon)}}\bigg(\sum_{j=-\infty}^{-1}b^{j\alpha(0)p(1+\varepsilon)}\|f\chi_j\|^{p(1+\varepsilon)}_{L^{q(\cdot)}(\mathbb{R}^{n})}\\
&\displaystyle\hspace{3mm}+\sum_{j=0}^\infty b^{j\alpha_\infty p(1+\varepsilon)}\|f\chi_j\|^{p(1+\varepsilon)}_{L^{q(\cdot)}(\mathbb{R}^{n})}\bigg)^{\frac{1}{p(1+\varepsilon)}}\\
&\displaystyle=\sup_{\varepsilon>0}\varepsilon^{\frac{\theta}{p(1+\varepsilon)}}\left(I_1+I_2\right)^{\frac{1}{p(1+\varepsilon)}},

\end{array}$$
where
$$I_1=\sum_{j=-\infty}^{-1}b^{j\alpha(0)p(1+\varepsilon)}\|f\chi_j\|^{p(1+\varepsilon)}_{L^{q(\cdot)}(\mathbb{R}^{n})},\,\,I_2=\sum_{j=0}^\infty b^{j\alpha_\infty p(1+\varepsilon)}\|f\chi_j\|^{p(1+\varepsilon)}_{L^{q(\cdot)}(\mathbb{R}^{n})}.$$
For $I_1$, take $\frac{1}{p(1+\varepsilon)}+\frac{1}{[p(1+\varepsilon)]'}=1$, by (3.2), $b>1, 0<\alpha(0), \alpha_\infty<\infty$ and the H\"{o}lder inequality, we have
$$\begin{array}{rl}
\displaystyle I_1&\displaystyle\lesssim \sum_{j=-\infty}^{-1}b^{j\alpha(0)p(1+\varepsilon)}\left(\sum_{k=j}^{-1}|\lambda_k|\|b_k\|_{L^{q(\cdot)}(\mathbb{R}^{n})}+\sum_{k=0}^\infty|\lambda_k|\|b_k\|_{L^{q(\cdot)}(\mathbb{R}^{n})}\right)^{p(1+\varepsilon)}\\
&\displaystyle\lesssim \sum_{j=-\infty}^{-1}\left(\sum_{k=j}^{-1}|\lambda_k|b^{(j-k)\alpha(0)}\right)^{p(1+\varepsilon)}+\sum_{j=-\infty}^{-1}\left(\sum_{k=0}^\infty|\lambda_k|b^{-k\alpha_\infty+j\alpha(0)}\right)^{p(1+\varepsilon)}\\
&\displaystyle\lesssim \sum_{j=-\infty}^{-1}\left(\sum_{k=j}^{-1}|\lambda_k|^{p(1+\varepsilon)}b^{(j-k)\alpha(0)p(1+\varepsilon)/2}\right)\left(\sum_{k=j}^{-1}b^{(j-k)\alpha(0)[p(1+\varepsilon)]'/2}\right)^{p(1+\varepsilon)/[p(1+\varepsilon)]'}\\
&\displaystyle\hspace{3mm}+\sum_{j=-\infty}^{-1}b^{j\alpha(0)p(1+\varepsilon)}\left(\sum_{k=0}^\infty|\lambda_k|^{p(1+\varepsilon)}b^{-k\alpha_\infty p(1+\varepsilon)/2}\right)\left(\sum_{k=0}^\infty b^{-k\alpha_\infty [p(1+\varepsilon)]'/2}\right)^{p(1+\varepsilon)/[p(1+\varepsilon)]'}\\

\end{array}$$
$$\begin{array}{rl}
&\displaystyle\lesssim\sum_{j=-\infty}^{-1}\left(\sum_{k=j}^{-1}|\lambda_k|^{p(1+\varepsilon)}b^{(j-k)\alpha(0)p(1+\varepsilon)/2}\right)+\sum_{k=0}^\infty|\lambda_k|^{p(1+\varepsilon)}b^{-k\alpha_\infty p(1+\varepsilon)/2}\\
&\displaystyle\lesssim\sum_{k=-\infty}^{-1}\sum_{j=-\infty}^k|\lambda_k|^{p(1+\varepsilon)}b^{(j-k)\alpha(0)p(1+\varepsilon)/2}+\sum_{k=0}^\infty|\lambda_k|^{p(1+\varepsilon)}\\
&\displaystyle\lesssim \sum_{k=-\infty}^{-1}|\lambda_k|^{p(1+\varepsilon)}+\sum_{k=0}^\infty|\lambda_k|^{p(1+\varepsilon)}\\
&\displaystyle\lesssim \sum_{k=-\infty}^\infty|\lambda_k|^{p(1+\varepsilon)}.

\end{array}$$
For $I_2$, take $\frac{1}{p(1+\varepsilon)}+\frac{1}{[p(1+\varepsilon)]'}=1$, by (3.2), $b>1, 0<\alpha_\infty<\infty$ and the H\"{o}lder inequality, we have
$$\begin{array}{rl}
\displaystyle I_2&\displaystyle\lesssim\sum_{j=0}^\infty b^{j\alpha_\infty p(1+\varepsilon)}\left(\sum_{k=j}^\infty|\lambda_k|\|b_k\|_{L^{q(\cdot)}(\mathbb{R}^{n})}\right)^{p(1+\varepsilon)}\\
&\displaystyle\lesssim \sum_{j=0}^\infty b^{j\alpha_\infty p(1+\varepsilon)}\left(\sum_{k=j}^\infty|\lambda_k|b^{-k\alpha_\infty}\right)^{p(1+\varepsilon)}\\
&\displaystyle\lesssim \sum_{j=0}^\infty b^{j\alpha_\infty p(1+\varepsilon)}\left(\sum_{k=j}^\infty|\lambda_k|^{p(1+\varepsilon)}b^{-k\alpha_\infty p(1+\varepsilon)/2}\right)\left(\sum_{k=j}^\infty b^{-k\alpha_\infty [p(1+\varepsilon)]'/2}\right)^{p(1+\varepsilon)/[p(1+\varepsilon)]'}\\
&\displaystyle\lesssim \sum_{j=0}^\infty\left(\sum_{k=j}^\infty|\lambda_k|^{p(1+\varepsilon)}b^{(j-k)\alpha_\infty p(1+\varepsilon)/2}\right)\left(\sum_{k=j}^\infty b^{(j-k)\alpha_\infty [p(1+\varepsilon)]'/2}\right)^{p(1+\varepsilon)/[p(1+\varepsilon)]'}\\
&\displaystyle\lesssim \sum_{k=0}^\infty|\lambda_k|^{p(1+\varepsilon)}\sum_{j=0}^kb^{(j-k)\alpha_\infty p(1+\varepsilon)/2}\\
&\displaystyle\lesssim \sum_{k=0}^\infty|\lambda_k|^{p(1+\varepsilon)}.

\end{array}$$
Thus we have
$$\|f\|_{\dot{K}^{\alpha(\cdot),p)}_{q(\cdot),\theta}(A;\mathbb{R}^{n})}\lesssim \sup_{\varepsilon>0}\varepsilon^{\frac{\theta}{p(1+\varepsilon)}}\left(\sum_{k=-\infty}^\infty|\lambda_k|^{p(1+\varepsilon)}\right)^{\frac{1}{p(1+\varepsilon)}}<\infty.$$
This leads to that
$f\in\dot{K}^{\alpha(\cdot),p)}_{q(\cdot),\theta}(A;\mathbb{R}^{n})$ and then
completes the proof of Theorem 3.1.

\noindent {\bf Remark 3.1}\quad From the proof of Theorem 3.1, it is
easy to see that if
$f\in\dot{K}^{\alpha(\cdot),p)}_{q(\cdot),\theta}(A;\mathbb{R}^{n})$ and
$\displaystyle f(x)=\sum_{k=-\infty}^\infty\lambda_kb_k(x)$ be a
central $(\alpha(\cdot), q(\cdot))$-block decomposition, then
$$\|f\|_{\dot{K}^{\alpha(\cdot),p)}_{q(\cdot),\theta}(A;\mathbb{R}^{n})}\approx\displaystyle\sup_{\varepsilon>0}\left(\varepsilon^\theta\sum_{k=-\infty}^\infty|\lambda_k|^{p(1+\varepsilon)}\right)^{\frac{1}{p(1+\varepsilon)}}.$$

By an argument similar to the proof of Theorem 3.1, we can obtain
the decomposition characterizations of the non-homogeneous anisotropic grand Herz
spaces with two variable exponents as follows.

\noindent{\bf Theorem 3.2}\quad Let $1\leq p<\infty, q(\cdot)\in \mathcal{P}(\mathbb{R}^{n}), \alpha(\cdot)\in L^\infty(\mathbb{R}^{n})\cap\mathcal{P}_0(\mathbb{R}^{n})\cap\mathcal{P}_\infty(\mathbb{R}^{n})$
and $0<\alpha_\infty<\infty$. The following two
statements are equivalent:

(i) $f\in K^{\alpha(\cdot),p)}_{q(\cdot),\theta}(A;\mathbb{R}^{n})$.

(ii) $f$ can be represented by
$$f(x)=\sum_{k=0}^\infty\lambda_kb_k(x),$$ where each $b_k$
is a central $(\alpha(\cdot), q(\cdot))$-block of restricted type
with support contained in $B_k$ and $\displaystyle\sup_{\varepsilon>0}\left(\varepsilon^\theta\sum_{k=0}^\infty|\lambda_k|^{p(1+\varepsilon)}\right)^{\frac{1}{p(1+\varepsilon)}}<\infty$.
\medskip

As applications of the decomposition theorems, let us come to
investigate the boundedness on the anisotropic grand Herz spaces with two variable
exponents for some sublinear operators.

\noindent{\bf Theorem 3.3}\quad Let $1\leq p<\infty, q(\cdot)\in \mathcal{B}(\mathbb{R}^{n}), \alpha(\cdot)\in L^\infty(\mathbb{R}^{n})\cap\mathcal{P}_0(\mathbb{R}^{n})\cap\mathcal{P}_\infty(\mathbb{R}^{n})$
and $0<\alpha(0), \alpha_\infty<\delta_2$. If a sublinear
operator $T$ satisfies $$|Tf(x)|\lesssim \int_{\mathbb{R}^{n}}\frac{|f(y)|}{\rho(x-y)}dy, \quad x\notin \mathrm{supp} f,\eqno(3.3)$$ for any
$f\in L^{q(\cdot)}(\mathbb{R}^{n})$ with a compact support and $T$ is bounded on
$L^{q(\cdot)}(\mathbb{R}^{n})$, then $T$ is bounded on
$\dot{K}^{\alpha(\cdot),p)}_{q(\cdot),\theta}(A;\mathbb{R}^{n})$ and
$K^{\alpha(\cdot),p)}_{q(\cdot),\theta}(A;\mathbb{R}^{n})$, respectively.

\noindent{\bf Proof}\quad It suffices to prove that $T$ is bounded on
$\dot{K}^{\alpha(\cdot),p)}_{q(\cdot),\theta}(A;\mathbb{R}^{n})$. The non-homogeneous case can be proved in the similar way. Suppose $f\in
\dot{K}^{\alpha(\cdot),p)}_{q(\cdot),\theta}(A;\mathbb{R}^{n})$. By Theorem 3.1,
$\displaystyle f(x)=\sum_{j=-\infty}^\infty\lambda_jb_j(x)$, where
each $b_j$ is a central $(\alpha(\cdot), q(\cdot))$-block with
support contained in $B_j$ and
$$\displaystyle\|f\|_{\dot{K}^{\alpha(\cdot),p)}_{q(\cdot),\theta}(A;\mathbb{R}^{n})}\approx\sup_{\varepsilon>0}\left(\varepsilon^\theta\sum_{j=-\infty}^\infty|\lambda_j|^{p(1+\varepsilon)}\right)^{\frac{1}{p(1+\varepsilon)}}.\eqno(3.4)$$
Therefore, we get
$$\begin{array}{rl}
&\displaystyle \sum_{k=-\infty}^\infty\|b^{k\alpha(\cdot)}(Tf)\chi_k\|^{p(1+\varepsilon)}_{L^{q(\cdot)}(\mathbb{R}^{n})}\\
&\displaystyle\lesssim\sum_{k=-\infty}^{-1}b^{k\alpha(0)p(1+\varepsilon)}\|(Tf)\chi_k\|^{p(1+\varepsilon)}_{L^{q(\cdot)}(\mathbb{R}^{n})}+\sum_{k=0}^\infty b^{k\alpha_\infty p(1+\varepsilon)}\|(Tf)\chi_k\|^{p(1+\varepsilon)}_{L^{q(\cdot)}(\mathbb{R}^{n})}\\
&\displaystyle\lesssim \sum_{k=-\infty}^{-1}b^{k\alpha(0) p(1+\varepsilon)}\bigg(\sum_{j=-\infty}^{k-w-1}|\lambda_j|\|(Tb_j)\chi_k\|_{L^{q(\cdot)}(\mathbb{R}^{n})}\bigg)^{p(1+\varepsilon)}\\
&\displaystyle\hspace{3mm}+\sum_{k=-\infty}^{-1}b^{k\alpha(0)
p(1+\varepsilon)}\bigg(\sum_{j=k-w}^{\infty}|\lambda_j|\|(Tb_j)\chi_k\|_{L^{q(\cdot)}(\mathbb{R}^{n})}\bigg)^{p(1+\varepsilon)}\\
&\displaystyle\hspace{3mm}+\sum_{k=0}^\infty b^{k\alpha_\infty
p(1+\varepsilon)}\bigg(\sum_{j=-\infty}^{k-w-1}|\lambda_j|\|(Tb_j)\chi_k\|_{L^{q(\cdot)}(\mathbb{R}^{n})}\bigg)^{p(1+\varepsilon)}\\
&\displaystyle\hspace{3mm}+\sum_{k=0}^\infty b^{k\alpha_\infty
p(1+\varepsilon)}\bigg(\sum_{j=k-w}^{\infty}|\lambda_j|\|(Tb_j)\chi_k\|_{L^{q(\cdot)}(\mathbb{R}^{n})}\bigg)^{p(1+\varepsilon)}\\
&\displaystyle=II_1+II_2+II_3+II_4.

\end{array}$$

Let us first estimate $II_1$. If $j\leq k-w-1$,\,$x\in C_k$ and $y\in B_j$, by (2.1) we have
$$\rho(x-y)\geq b^{-w}\rho(x)-\rho(y)\geq b^{-w}\rho(x)-b^{-w-1}\rho(x)=b^{-w}(1-1/b)\rho(x).\eqno(3.5)$$
Therefore by (3.3) and the generalized H\"{o}lder
inequality, we get
$$\begin{array}{rl}
\displaystyle |Tb_j(x)|&\displaystyle\leq C\rho(x)^{-1}\int_{B_j}|b_j(y)|dy\\
&\displaystyle\leq Cb^{-k}\|b_j\|_{L^{q(\cdot)}(\mathbb{R}^{n})}\|\chi_{B_j}\|_{L^{q'(\cdot)}(\mathbb{R}^{n})}.

\end{array}$$
So by Lemma 1.2 and Lemma 1.3, we have
$$\begin{array}{rl}
\displaystyle \|(Tb_j)\chi_k\|_{L^{q(\cdot)}(\mathbb{R}^{n})}&\displaystyle\lesssim b^{-k}\|b_j\|_{L^{q(\cdot)}(\mathbb{R}^{n})}\|\chi_{B_j}\|_{L^{q'(\cdot)}(\mathbb{R}^{n})}\|\chi_{B_k}\|_{L^{q(\cdot)}(\mathbb{R}^{n})}\\
&\displaystyle\lesssim b^{-k}\|b_j\|_{L^{q(\cdot)}(\mathbb{R}^{n})}\big(|B_k|\|\chi_{B_k}\|^{-1}_{L^{q'(\cdot)}(\mathbb{R}^{n})}\big)\|\chi_{B_j}\|_{L^{q'(\cdot)}(\mathbb{R}^{n})}\\
&\displaystyle\lesssim\|b_j\|_{L^{q(\cdot)}(\mathbb{R}^{n})}\frac{\|\chi_{B_j}\|_{L^{q'(\cdot)}(\mathbb{R}^{n})}}{\|\chi_{B_k}\|_{L^{q'(\cdot)}(\mathbb{R}^{n})}}\\
&\displaystyle\lesssim b^{\delta_2(j-k)}\|b_j\|_{L^{q(\cdot)}(\mathbb{R}^{n})}.

\end{array}\eqno(3.6)$$
Since $b>1, 0<\alpha(0)<\delta_2$, by
(3.6) and the H\"{o}lder inequality, we have
$$\begin{array}{rl}
\displaystyle II_1&\displaystyle\lesssim\sum_{k=-\infty}^{-1}b^{k\alpha(0) p(1+\varepsilon)}\bigg(\sum_{j=-\infty}^{k-w-1}|\lambda_j|b^{\delta_2(j-k)-j\alpha(0)}\bigg)^{p(1+\varepsilon)}\\
&\displaystyle\lesssim\sum_{k=-\infty}^{-1}\bigg(\sum_{j=-\infty}^{k-w-1}|\lambda_j|^{p(1+\varepsilon)}b^{(j-k)[\delta_2-\alpha(0)]p(1+\varepsilon)/2}\bigg)
\bigg(\sum_{j=-\infty}^{k-w-1}b^{(j-k)[\delta_2-\alpha(0)][p(1+\varepsilon)]'/2}\bigg)^{p(1+\varepsilon)/[p(1+\varepsilon)]'}\\
&\displaystyle\lesssim\sum_{k=-\infty}^{-1}\bigg(\sum_{j=-\infty}^{k-w-1}|\lambda_j|^{p(1+\varepsilon)}b^{(j-k)[\delta_2-\alpha(0)]p(1+\varepsilon)/2}\bigg)\\
&\displaystyle\lesssim\sum_{j=-\infty}^{-w-2}|\lambda_j|^{p(1+\varepsilon)}\sum_{k=j+w+1}^{-1}b^{(j-k)[\delta_2-\alpha(0)]p(1+\varepsilon)/2}\\
&\displaystyle\lesssim\sum_{j=-\infty}^{-w-2}|\lambda_j|^{p(1+\varepsilon)}.

\end{array}\eqno(3.7)$$

Let us now estimate $II_2$. Since $b>1, 0<\alpha(0), \alpha_\infty<\delta_2$, by
$L^{q(\cdot)}(\mathbb{R}^{n})$ boundedness of $T$ and the H\"{o}lder
inequality, we have
$$\begin{array}{rl}
\displaystyle II_2&\displaystyle\lesssim\sum_{k=-\infty}^{-1}b^{k\alpha(0) p(1+\varepsilon)}\bigg(\sum_{j=k-w}^{\infty}|\lambda_j|\|b_j\|_{L^{q(\cdot)}(\mathbb{R}^{n})}\bigg)^{p(1+\varepsilon)}\\
&\displaystyle\lesssim\sum_{k=-\infty}^{-1}\bigg(\sum_{j=k-w}^{-1}|\lambda_j|b^{(k-j)\alpha(0)}\bigg)^{p(1+\varepsilon)}+\sum_{k=-\infty}^{-1}b^{k\alpha(0) p(1+\varepsilon)}\bigg(\sum_{j=0}^{\infty}|\lambda_j|b^{-j\alpha_\infty}\bigg)^{p(1+\varepsilon)}\\
&\displaystyle\lesssim\sum_{k=-\infty}^{-1}\bigg(\sum_{j=k-w}^{-1}|\lambda_j|^{p(1+\varepsilon)}b^{(k-j)\alpha(0)p(1+\varepsilon)/2}\bigg)\bigg(\sum_{j=k-w}^{-1}b^{(k-j)\alpha(0)[p(1+\varepsilon)]'/2}\bigg)^{p(1+\varepsilon)/[p(1+\varepsilon)]'}\\
&\displaystyle\hspace{3mm}+\bigg(\sum_{j=0}^{\infty}|\lambda_j|^{p(1+\varepsilon)}b^{-j\alpha_\infty p(1+\varepsilon)/2}\bigg)\bigg(\sum_{j=0}^{\infty}b^{-j\alpha_\infty [p(1+\varepsilon)]'/2}\bigg)^{p(1+\varepsilon)/[p(1+\varepsilon)]'}\\
&\displaystyle\lesssim\sum_{j=-\infty}^{-1}|\lambda_j|^{p(1+\varepsilon)}\sum_{k=-\infty}^{j+w}b^{(k-j)\alpha(0)p(1+\varepsilon)/2}+\sum_{j=0}^{\infty}|\lambda_j|^{p(1+\varepsilon)}\\

\end{array}$$
$$\begin{array}{rl}
&\displaystyle\lesssim\sum_{j=-\infty}^\infty|\lambda_j|^{p(1+\varepsilon)}.

\end{array}\eqno(3.8)$$

For $II_3$. Since $0<\alpha(0), \alpha_\infty<\delta_2$, by
(3.6) and the H\"{o}lder inequality, we have
$$\begin{array}{rl}
\displaystyle II_3&\displaystyle\lesssim\sum_{k=0}^\infty b^{k\alpha_\infty p(1+\varepsilon)}\bigg(\sum_{j=-\infty}^{k-w-1}|\lambda_j|b^{\delta_2(j-k)}\|b_j\|_{L^{q(\cdot)}(\mathbb{R}^{n})}\bigg)^{p(1+\varepsilon)}\\
&\displaystyle\lesssim\sum_{k=0}^\infty b^{k\alpha_\infty p(1+\varepsilon)}\bigg(\sum_{j=-\infty}^{-1}|\lambda_j|b^{\delta_2(j-k)-j\alpha(0)}\bigg)^{p(1+\varepsilon)}\\
&\displaystyle\hspace{3mm}+\sum_{k=0}^\infty b^{k\alpha_\infty p(1+\varepsilon)}\bigg(\sum_{j=0}^{k-w-1}|\lambda_j|b^{\delta_2(j-k)-j\alpha_\infty}\bigg)^{p(1+\varepsilon)}\\
&\displaystyle\lesssim\sum_{k=0}^\infty b^{k(\alpha_\infty-\delta_2) p(1+\varepsilon)}\bigg(\sum_{j=-\infty}^{-1}|\lambda_j|b^{j[\delta_2-\alpha(0)]}\bigg)^{p(1+\varepsilon)}\\
&\displaystyle\hspace{3mm}+\sum_{k=0}^\infty\bigg(\sum_{j=0}^{k-w-1}|\lambda_j|b^{(j-k)(\delta_2-\alpha_\infty)}\bigg)^{p(1+\varepsilon)}\\
&\displaystyle\lesssim\bigg(\sum_{j=-\infty}^{-1}|\lambda_j|^{p(1+\varepsilon)}b^{j[\delta_2-\alpha(0)]p(1+\varepsilon)/2}\bigg)\bigg(\sum_{j=-\infty}^{-1}b^{j[\delta_2-\alpha(0)][p(1+\varepsilon)]'/2}\bigg)^{p(1+\varepsilon)/[p(1+\varepsilon)]'}\\
&\displaystyle\hspace{3mm}+\sum_{k=0}^\infty\bigg(\sum_{j=0}^{k-w-1}|\lambda_j|^{p(1+\varepsilon)}b^{(j-k)(\delta_2-\alpha_\infty)p(1+\varepsilon)/2}\bigg)\bigg(\sum_{j=0}^{k-w-1}b^{(j-k)(\delta_2-\alpha_\infty)[p(1+\varepsilon)]'/2}\bigg)^{p(1+\varepsilon)/[p(1+\varepsilon)]'}\\
&\displaystyle\lesssim\sum_{j=-\infty}^{-1}|\lambda_j|^{p(1+\varepsilon)}b^{j[\delta_2-\alpha(0)]p(1+\varepsilon)/2}+\sum_{k=0}^\infty\sum_{j=0}^{k-w-1}|\lambda_j|^{p(1+\varepsilon)}b^{(j-k)(\delta_2-\alpha_\infty)p(1+\varepsilon)/2}\\
&\displaystyle\lesssim\sum_{j=-\infty}^{-1}|\lambda_j|^{p(1+\varepsilon)}+\sum_{j=0}^\infty|\lambda_j|^{p(1+\varepsilon)}\sum_{k=j+w+1}^\infty b^{(j-k)(\delta_2-\alpha_\infty)p(1+\varepsilon)/2}\\
&\displaystyle\lesssim\sum_{j=-\infty}^\infty|\lambda_j|^{p(1+\varepsilon)}.

\end{array}\eqno(3.9)$$

Let us now estimate $II_4$. By
$L^{q(\cdot)}(\mathbb{R}^{n})$ boundedness of $T$ and the H\"{o}lder
inequality, we have
$$\begin{array}{rl}
\displaystyle II_4&\displaystyle=\sum_{k=0}^\infty b^{k\alpha_\infty p(1+\varepsilon)}\bigg(\sum_{j=k-w}^{\infty}|\lambda_j|\|(Tb_j)\chi_k\|_{L^{q(\cdot)}(\mathbb{R}^{n})}\bigg)^{p(1+\varepsilon)}\\
&\displaystyle\lesssim\sum_{k=0}^\infty b^{k\alpha_\infty p(1+\varepsilon)}\bigg(\sum_{j=k-w}^{\infty}|\lambda_j|\|b_j\|_{L^{q(\cdot)}(\mathbb{R}^{n})}\bigg)^{p(1+\varepsilon)}\\
&\displaystyle\lesssim\sum_{k=0}^\infty b^{k\alpha_\infty p(1+\varepsilon)}\bigg(\sum_{j=k-w}^\infty|\lambda_j|b^{-j\alpha_\infty}\bigg)^{p(1+\varepsilon)}\\
&\displaystyle\lesssim\sum_{k=0}^\infty b^{k\alpha_\infty p(1+\varepsilon)}\bigg(\sum_{j=k-w}^\infty|\lambda_j|^{p(1+\varepsilon)}b^{-j\alpha_\infty p(1+\varepsilon)/2}\bigg)\bigg(\sum_{j=k-w}^\infty b^{-j\alpha_\infty [p(1+\varepsilon)]'/2}\bigg)^{p(1+\varepsilon)/[p(1+\varepsilon)]'}\\
&\displaystyle\lesssim\sum_{k=0}^\infty b^{k\alpha_\infty p(1+\varepsilon)/2}\bigg(\sum_{j=k-w}^\infty|\lambda_j|^{p(1+\varepsilon)}b^{-j\alpha_\infty p(1+\varepsilon)/2}\bigg)\\
&\displaystyle\lesssim\sum_{j=-w}^\infty|\lambda_j|^{p(1+\varepsilon)}\sum_{k=0}^{j+w}b^{(k-j)\alpha_\infty p(1+\varepsilon)/2}\\

\end{array}$$
$$\begin{array}{rl}
&\displaystyle\lesssim\sum_{j=-w}^\infty|\lambda_j|^{p(1+\varepsilon)}.

\end{array}\eqno(3.10)$$

Combining (3.7)-(3.10), by (3.4) we have
$$\begin{array}{rl}
\displaystyle\|Tf\|_{\dot{K}^{\alpha(\cdot),p)}_{q(\cdot),\theta}(A;\mathbb{R}^{n})}&\displaystyle=\sup_{\varepsilon>0}\left(\varepsilon^{\theta}\sum_{k=-\infty}^\infty\|b^{k\alpha(\cdot)}(Tf)\chi_k\|^{p(1+\varepsilon)}_{L^{q(\cdot)}(\mathbb{R}^{n})}\right)^{\frac{1}{p(1+\varepsilon)}}\\
&\displaystyle\lesssim\sup_{\varepsilon>0}\left(\varepsilon^{\theta}\sum_{j=-\infty}^\infty|\lambda_j|^{p(1+\varepsilon)}\right)^{\frac{1}{p(1+\varepsilon)}}\\
&\displaystyle\lesssim\|f\|_{\dot{K}^{\alpha(\cdot),p)}_{q(\cdot),\theta}(A;\mathbb{R}^{n})}.
\end{array}$$

Thus, the proof of Theorem 3.3 is completed.

\noindent {\bf Remark 3.1}\quad From the proof of Theorem 3.3, it is
easy to see that the size condition (3.3) can be replaced by
$$|Tf(x)|\leq C\frac{\|f\|_{L^1}}{\rho(x)},\quad \mathrm{if}\,\,\inf_{y\in \mathrm{supp}f}\rho(x-y)\geq b^{-w}(1-1/b)\rho(x),\eqno(3.11)$$
for all $f\in L^{q(\cdot)}(\mathbb{R}^n)$ with compact support.

\section{Anisotropic grand Herz-Morrey spaces with variable
exponents\label{s4}}$\indent$

In this section, we first introduce the definition of anisotropic grand Herz-Morrey spaces with two variable
exponents.

\noindent {\bf Definition 4.1}\quad Let $\alpha(\cdot):
\mathbb{R}^{n}\rightarrow\mathbb{R}$ with $\alpha(\cdot)\in L^\infty(\mathbb{R}^{n})$, $1\leq p<\infty$, $q(\cdot)\in
\mathcal{P}(\mathbb{R}^{n})$, $0\leq \lambda<\infty$ and $\theta>0$. The homogeneous anisotropic grand Herz-Morrey space
$M\dot{K}^{\alpha(\cdot),\lambda,\theta}_{q(\cdot),p)}(A;\mathbb{R}^{n})$ associated with the dilation $A$ is defined by
$$M\dot{K}^{\alpha(\cdot),\lambda,\theta}_{q(\cdot),p)}(A;\mathbb{R}^{n})=\{f\in L_{\rm
loc}^{q(\cdot)}(\mathbb{R}^{n}\setminus \{0\}):
\|f\|_{M\dot{K}^{\alpha(\cdot),\lambda,\theta}_{q(\cdot),p)}(A;\mathbb{R}^{n})}<\infty\},$$
where
$$\|f\|_{M\dot{K}^{\alpha(\cdot),\lambda,\theta}_{q(\cdot),p)}(A;\mathbb{R}^{n})}=\sup_{\varepsilon>0}\sup_{L\in\mathbb{Z}}b^{-L\lambda}\left\{\varepsilon^\theta\sum_{k=-\infty}^L\|b^{k\alpha(\cdot)
}f\chi_k\|^{p(1+\varepsilon)}_{L^{q(\cdot)}(\mathbb{R}^{n})}\right\}^{\frac{1}{p(1+\varepsilon)}}.$$ The
non-homogeneous anisotropic grand Herz-Morrey space $MK^{\alpha(\cdot),\lambda,\theta}_{q(\cdot),p)}(A;\mathbb{R}^{n})$ associated with the dilation $A$ is defined by $$MK^{\alpha(\cdot),\lambda,\theta}_{q(\cdot),p)}(A;\mathbb{R}^{n})=\{f\in
L_{\rm loc}^{q(\cdot)}(\mathbb{R}^{n}):
\|f\|_{MK^{\alpha(\cdot),\lambda,\theta}_{q(\cdot),p)}(A;\mathbb{R}^{n})}<\infty\},$$ where
$$\|f\|_{MK^{\alpha(\cdot),\lambda,\theta}_{q(\cdot),p)}(A;\mathbb{R}^{n})}=\sup_{\varepsilon>0}\sup_{L\in\mathbb{Z}}b^{-L\lambda}\left\{\varepsilon^\theta\sum_{k=0}^L\|b^{k\alpha(\cdot)}f\tilde{\chi}_k\|^{p(1+\varepsilon)}_{L^{q(\cdot)}(\mathbb{R}^{n})}\right\}^{\frac{1}{p(1+\varepsilon)}}.$$

\noindent {\bf Remark 4.1}\quad If $\lambda=0$, then $$M\dot{K}^{\alpha(\cdot),0,\theta}_{q(\cdot),p)}(A;\mathbb{R}^{n})=\dot{K}^{\alpha(\cdot),p)}_{q(\cdot),\theta}(A;\mathbb{R}^{n})$$ and $$MK^{\alpha(\cdot),0,\theta}_{q(\cdot),p)}(A;\mathbb{R}^{n})=K^{\alpha(\cdot),p)}_{q(\cdot),\theta}(A;\mathbb{R}^{n}),$$
where $\dot{K}^{\alpha(\cdot),p)}_{q(\cdot),\theta}(A;\mathbb{R}^{n})$ and $K^{\alpha(\cdot),p)}_{q(\cdot),\theta}(A;\mathbb{R}^{n})$ are the anisotropic Herz spaces associated with the dilation $A$.

From Lemma 3.1 we can easily get the following lemma.

\noindent{\bf Lemma 4.1}\quad Let $1\leq p<\infty$, $q(\cdot)\in
\mathcal{P}(\mathbb{R}^{n})$, $\theta>0$ and $\alpha(\cdot)\in L^\infty(\mathbb{R}^{n})\cap\mathcal{P}_\infty(\mathbb{R}^{n})$. Then
$$\begin{array}{rl}
\displaystyle \|f\|_{M\dot{K}^{\alpha(\cdot),\lambda,\theta}_{q(\cdot),p)}(A;\mathbb{R}^{n})}&\displaystyle\approx\max\bigg\{\sup_{\varepsilon>0}\sup_{L\leq0,L\in\mathbb{Z}}b^{-L\lambda}\bigg(\varepsilon^\theta\sum_{k=-\infty}^L b^{k\alpha(0)p(1+\varepsilon)}\|f\chi_k\|^{p(1+\varepsilon)}_{L^{q(\cdot)}(\mathbb{R}^{n})}\bigg)^{\frac{1}{p(1+\varepsilon)}},\\
&\displaystyle\hspace{12mm}\sup_{\varepsilon>0}\sup_{L>0,L\in\mathbb{Z}}\bigg[b^{-L\lambda}\bigg(\varepsilon^\theta\sum_{k=-\infty}^{-1} b^{k\alpha(0)p(1+\varepsilon)}\|f\chi_k\|^{p(1+\varepsilon)}_{L^{q(\cdot)}(\mathbb{R}^{n})}\bigg)^{\frac{1}{p(1+\varepsilon)}}\\
&\displaystyle\hspace{12mm}+b^{-L\lambda}\bigg(\varepsilon^\theta\sum_{k=0}^L b^{k\alpha_\infty p(1+\varepsilon)}\|f\chi_k\|^{p(1+\varepsilon)}_{L^{q(\cdot)}(\mathbb{R}^{n})}\bigg)^{\frac{1}{p(1+\varepsilon)}}\bigg]\bigg\}.

\end{array}$$

Next we will consider some properties of
$M\dot{K}^{\alpha(\cdot),\lambda,\theta}_{q(\cdot),p)}(A;\mathbb{R}^{n})$. There are similar properties for $MK^{\alpha(\cdot),\lambda,\theta}_{q(\cdot),p)}(A;\mathbb{R}^{n})$.

\noindent{\bf Theorem 4.1}\quad Suppose $\alpha(\cdot),\alpha_1(\cdot),\alpha_2(\cdot)\in L^\infty(\mathbb{R}^{n})\cap\mathcal{P}_0(\mathbb{R}^{n})\cap\mathcal{P}_\infty(\mathbb{R}^{n})$, $q(\cdot),q_1(\cdot),q_2(\cdot)\in
\mathcal{P}(\mathbb{R}^{n})$, $1<p,p_1,p_2<\infty$, and $0\leq \lambda, \lambda_1, \lambda_2<\infty$ such that $\alpha(\cdot)=\alpha_1(\cdot)+\alpha_2(\cdot)$, $1/q(\cdot)=1/q_1(\cdot)+1/q_2(\cdot)$, $1/p=1/p_1+1/p_2$ and $\lambda=\lambda_1+\lambda_2$. If $f\in M\dot{K}^{\alpha_1(\cdot),\lambda_1,\theta}_{q_1(\cdot),p_1)}(A;\mathbb{R}^{n})$
and $g\in M\dot{K}^{\alpha_2(\cdot),\lambda_2,\theta}_{q_2(\cdot),p_2)}(A;\mathbb{R}^{n})$, then $fg\in M\dot{K}^{\alpha(\cdot),\lambda,\theta}_{q(\cdot),p)}(A;\mathbb{R}^{n})$ and
$$\|fg\|_{M\dot{K}^{\alpha(\cdot),\lambda,\theta}_{q(\cdot),p)}(A;\mathbb{R}^{n})} \leq \|f\|_{M\dot{K}^{\alpha_1(\cdot),\lambda_1,\theta}_{q_1(\cdot),p_1)}(A;\mathbb{R}^{n})}\|g\|_{M\dot{K}^{\alpha_2(\cdot),\lambda_2,\theta}_{q_2(\cdot),p_2)}(A;\mathbb{R}^{n})}.$$

\noindent{\bf Proof}\quad By Lemma 2.4 and the H\"{o}lder inequality we have
$$\begin{array}{rl}
\displaystyle \|fg\|_{M\dot{K}^{\alpha(\cdot),\lambda,\theta}_{q(\cdot),p)}(A;\mathbb{R}^{n})}&\displaystyle=\sup_{\varepsilon>0}\sup_{L\in\mathbb{Z}}b^{-L\lambda}\left(\varepsilon^\theta\sum_{k=-\infty}^L\|b^{k\alpha(\cdot)
}fg\chi_k\|^{p(1+\varepsilon)}_{L^{q(\cdot)}(\mathbb{R}^{n})}\right)^{\frac{1}{p(1+\varepsilon)}}\\
&\displaystyle\leq\sup_{\varepsilon>0}\sup_{L\in\mathbb{Z}}b^{-L\lambda_1}b^{-L\lambda_2}\left(\varepsilon^\theta\sum_{k=-\infty}^L\|b^{k\alpha_1(\cdot)
}f\chi_k\|^{p(1+\varepsilon)}_{L^{q_1(\cdot)}(\mathbb{R}^{n})}\|b^{k\alpha_2(\cdot)
}g\chi_k\|^{p(1+\varepsilon)}_{L^{q_2(\cdot)}(\mathbb{R}^{n})}\right)^{\frac{1}{p(1+\varepsilon)}}\\
&\displaystyle\leq\sup_{\varepsilon>0}\sup_{L\in\mathbb{Z}}b^{-L\lambda_1}\left(\varepsilon^\theta\sum_{k=-\infty}^L\|b^{k\alpha_1(\cdot)
}f\chi_k\|^{p_1(1+\varepsilon)}_{L^{q_1(\cdot)}(\mathbb{R}^{n})}\right)^{\frac{1}{p_1(1+\varepsilon)}}\\
&\displaystyle\hspace{3mm}\times \sup_{\varepsilon>0}\sup_{L\in\mathbb{Z}}b^{-L\lambda_2}\left(\varepsilon^\theta\sum_{k=-\infty}^L\|b^{k\alpha_2(\cdot)
}g\chi_k\|^{p_2(1+\varepsilon)}_{L^{q_2(\cdot)}(\mathbb{R}^{n})}\right)^{\frac{1}{p_2(1+\varepsilon)}}\\
&\displaystyle=\|f\|_{M\dot{K}^{\alpha_1(\cdot),\lambda_1,\theta}_{q_1(\cdot),p_1)}(A;\mathbb{R}^{n})}\|g\|_{M\dot{K}^{\alpha_2(\cdot),\lambda_2,\theta}_{q_2(\cdot),p_2)}(A;\mathbb{R}^{n})}.

\end{array}$$

So we complete the proof of Theorem 4.1.

From Theorem 4.1 we can get the following corollary.

\noindent{\bf Corollary 4.1}\quad Suppose $\alpha(\cdot),\alpha_i(\cdot)\in L^\infty(\mathbb{R}^{n})\cap\mathcal{P}_0(\mathbb{R}^{n})\cap\mathcal{P}_\infty(\mathbb{R}^{n})$, $q(\cdot),q_i(\cdot)\in
\mathcal{P}(\mathbb{R}^{n})$, $1<p,p_i<\infty$, and $0\leq \lambda, \lambda_i<\infty(1\leq i\leq m)$ such that $\displaystyle\alpha(\cdot)=\sum_{i=1}^m\alpha_i(\cdot)$, $\displaystyle1/q(\cdot)=\sum_{i=1}^m1/q_i(\cdot)$, $\displaystyle1/p=\sum_{i=1}^m1/p_i$ and $\displaystyle\lambda=\sum_{i=1}^m\lambda_i$. If $f_i\in M\dot{K}^{\alpha_i(\cdot),\lambda_i,\theta}_{q_i(\cdot),p_i)}(A;\mathbb{R}^{n})$, then $\displaystyle\prod_{i=1}^mf_i\in M\dot{K}^{\alpha(\cdot),\lambda,\theta}_{q(\cdot),p)}(A;\mathbb{R}^{n})$ and
$$\left\|\prod_{i=1}^mf_i\right\|_{M\dot{K}^{\alpha(\cdot),\lambda,\theta}_{q(\cdot),p)}(A;\mathbb{R}^{n})} \leq \prod_{i=1}^m\|f_i\|_{M\dot{K}^{\alpha_i(\cdot),\lambda_i,\theta}_{q_i(\cdot),p_i)}(A;\mathbb{R}^{n})}.$$

\noindent{\bf Theorem 4.2}\quad Suppose $\alpha(\cdot)\in L^\infty(\mathbb{R}^{n})\cap\mathcal{P}_0(\mathbb{R}^{n})\cap\mathcal{P}_\infty(\mathbb{R}^{n})$, $q(\cdot)\in
\mathcal{P}(\mathbb{R}^{n})$, $1<p<\infty$, and $0\leq \lambda<\infty$. If $f,g\in M\dot{K}^{\alpha(\cdot),\lambda,\theta}_{q(\cdot),p)}(A;\mathbb{R}^{n})$, then $f+g\in M\dot{K}^{\alpha(\cdot),\lambda,\theta}_{q(\cdot),p)}(A;\mathbb{R}^{n})$ and
$$\|f+g\|_{M\dot{K}^{\alpha(\cdot),\lambda,\theta}_{q(\cdot),p)}(A;\mathbb{R}^{n})} \leq \|f\|_{M\dot{K}^{\alpha(\cdot),\lambda,\theta}_{q(\cdot),p)}(A;\mathbb{R}^{n})}+\|g\|_{M\dot{K}^{\alpha(\cdot),\lambda,\theta}_{q(\cdot),p)}(A;\mathbb{R}^{n})}.$$

\noindent{\bf Proof}\quad By the Minkowski inequality we have
$$\begin{array}{rl}
&\displaystyle \|f+g\|_{M\dot{K}^{\alpha(\cdot),\lambda,\theta}_{q(\cdot),p)}(A;\mathbb{R}^{n})}\\
&\displaystyle=\sup_{\varepsilon>0}\sup_{L\in\mathbb{Z}}b^{-L\lambda}\left(\varepsilon^\theta\sum_{k=-\infty}^L\|b^{k\alpha(\cdot)
}(f+g)\chi_k\|^{p(1+\varepsilon)}_{L^{q(\cdot)}(\mathbb{R}^{n})}\right)^{\frac{1}{p(1+\varepsilon)}}\\
&\displaystyle\leq\sup_{\varepsilon>0}\sup_{L\in\mathbb{Z}}b^{-L\lambda}\left[\varepsilon^\theta\sum_{k=-\infty}^L\left(\|b^{k\alpha(\cdot)
}f\chi_k\|^{p(1+\varepsilon)}_{L^{q(\cdot)}(\mathbb{R}^{n})}+\|b^{k\alpha(\cdot)
}g\chi_k\|^{p(1+\varepsilon)}_{L^{q(\cdot)}(\mathbb{R}^{n})}\right)\right]^{\frac{1}{p(1+\varepsilon)}}\\
&\displaystyle\leq\sup_{\varepsilon>0}\sup_{L\in\mathbb{Z}}b^{-L\lambda}\left[\varepsilon^\theta\left(\sum_{k=-\infty}^L\|b^{k\alpha(\cdot)
}f\chi_k\|^{p(1+\varepsilon)}_{L^{q(\cdot)}(\mathbb{R}^{n})}\right)^{\frac{1}{p(1+\varepsilon)}}+\varepsilon^\theta\left(\sum_{k=-\infty}^L\|b^{k\alpha(\cdot)
}g\chi_k\|^{p(1+\varepsilon)}_{L^{q(\cdot)}(\mathbb{R}^{n})}\right)^{\frac{1}{p(1+\varepsilon)}}\right]\\
&\displaystyle\leq\|f\|_{M\dot{K}^{\alpha(\cdot),\lambda,\theta}_{q(\cdot),p)}(A;\mathbb{R}^{n})}+\|g\|_{M\dot{K}^{\alpha(\cdot),\lambda,\theta}_{q(\cdot),p)}(A;\mathbb{R}^{n})}.

\end{array}$$

So we obtain Theorem 4.2.

From Theorem 4.2 we can get the following corollary.

\noindent{\bf Corollary 4.2}\quad Suppose $\alpha(\cdot)\in L^\infty(\mathbb{R}^{n})\cap\mathcal{P}_0(\mathbb{R}^{n})\cap\mathcal{P}_\infty(\mathbb{R}^{n})$, $q(\cdot)\in
\mathcal{P}(\mathbb{R}^{n})$, $1<p<\infty$, and $0\leq \lambda<\infty$. If $f_i\in M\dot{K}^{\alpha(\cdot),\lambda,\theta}_{q(\cdot),p)}(A;\mathbb{R}^{n}),1\leq i\leq m$, then $\displaystyle\sum_{i=1}^mf_i\in M\dot{K}^{\alpha(\cdot),\lambda,\theta}_{q(\cdot),p)}(A;\mathbb{R}^{n})$ and
$$\left\|\sum_{i=1}^mf_i\right\|_{M\dot{K}^{\alpha(\cdot),\lambda,\theta}_{q(\cdot),p)}(A;\mathbb{R}^{n})} \leq \sum_{i=1}^m\|f_i\|_{M\dot{K}^{\alpha(\cdot),\lambda,\theta}_{q(\cdot),p)}(A;\mathbb{R}^{n})}.$$

At the end of this section, we will
investigate the boundedness on the anisotropic Herz-Morrey spaces with two variable
exponents for some sublinear operators.

\noindent{\bf Theorem 4.3}\quad Let $1\leq p<\infty, q(\cdot)\in \mathcal{B}(\mathbb{R}^{n}), \alpha(\cdot)\in L^\infty(\mathbb{R}^{n})\cap\mathcal{P}_0(\mathbb{R}^{n})\cap\mathcal{P}_\infty(\mathbb{R}^{n})$,
$0<2\lambda<\alpha(0), \alpha_\infty<\delta_2$. If a sublinear
operator $T$ satisfies the condition (3.3) for any
$f\in L^{q(\cdot)}(\mathbb{R}^{n})$ with a compact support and $T$ is bounded on
$L^{q(\cdot)}(\mathbb{R}^{n})$, then $T$ is bounded on
$M\dot{K}^{\alpha(\cdot),\lambda,\theta}_{q(\cdot),p)}(A;\mathbb{R}^{n})$ and
$MK^{\alpha(\cdot),\lambda,\theta}_{q(\cdot),p)}(A;\mathbb{R}^{n})$, respectively.

\noindent{\bf Proof}\quad It suffices to prove that $T$ is bounded on
$M\dot{K}^{\alpha(\cdot),\lambda,\theta}_{q(\cdot),p)}(A;\mathbb{R}^{n})$. The non-homogeneous case can be proved in the similar way. Suppose $f\in M\dot{K}^{\alpha(\cdot),\lambda,\theta}_{q(\cdot),p)}(A;\mathbb{R}^{n})$. By Lemma 4.1 we get
$$\begin{array}{rl}
\displaystyle \|Tf\|_{M\dot{K}^{\alpha(\cdot),\lambda,\theta}_{q(\cdot),p)}(A;\mathbb{R}^{n})}&\displaystyle\approx\max\bigg\{\sup_{\varepsilon>0}\sup_{L\leq0,L\in\mathbb{Z}}b^{-L\lambda}\bigg(\varepsilon^\theta\sum_{k=-\infty}^L b^{k\alpha(0)p(1+\varepsilon)}\|(Tf)\chi_k\|^{p(1+\varepsilon)}_{L^{q(\cdot)}(\mathbb{R}^{n})}\bigg)^{\frac{1}{p(1+\varepsilon)}},\\
&\displaystyle\hspace{3mm}\sup_{\varepsilon>0}\sup_{L>0,L\in\mathbb{Z}}\bigg[b^{-L\lambda}\bigg(\varepsilon^\theta\sum_{k=-\infty}^{-1} b^{k\alpha(0)p(1+\varepsilon)}\|(Tf)\chi_k\|^{p(1+\varepsilon)}_{L^{q(\cdot)}(\mathbb{R}^{n})}\bigg)^{\frac{1}{p(1+\varepsilon)}}\\
&\displaystyle\hspace{3mm}+b^{-L\lambda}\bigg(\varepsilon^\theta\sum_{k=0}^L b^{k\alpha_\infty p(1+\varepsilon)}\|(Tf)\chi_k\|^{p(1+\varepsilon)}_{L^{q(\cdot)}(\mathbb{R}^{n})}\bigg)^{\frac{1}{p(1+\varepsilon)}}\bigg]\bigg\}\\
&\displaystyle=\max\{J_1, J_2\}.

\end{array}$$
It suffices to prove that $J_1$ is bounded on
$M\dot{K}^{\alpha(\cdot),\lambda,\theta}_{q(\cdot),p)}(A;\mathbb{R}^{n})$, as the estimate of $J_2$ is essentially similar to that of $J_1$. Denote $f_j=f\chi_j$ for each $j\in \mathbb{Z}$, then we have
$\displaystyle f(x)=\sum_{j=-\infty}^\infty f_j(x)$. It is easy to see that

$$\begin{array}{rl}
\displaystyle J_1&\displaystyle\lesssim\sup_{\varepsilon>0}\sup_{L\leq0,L\in\mathbb{Z}}b^{-L\lambda}\left[\varepsilon^\theta\sum_{k=-\infty}^L b^{k\alpha(0)p(1+\varepsilon)}\left(\sum_{j=-\infty}^{k-w-1}\|(Tf_j)\chi_k\|_{L^{q(\cdot)}(\mathbb{R}^{n})}\right)^{p(1+\varepsilon)}\right]^{\frac{1}{p(1+\varepsilon)}}\\
&\displaystyle\hspace{3mm}+\sup_{\varepsilon>0}\sup_{L\leq0,L\in\mathbb{Z}}b^{-L\lambda}\left[\varepsilon^\theta\sum_{k=-\infty}^L b^{k\alpha(0)p(1+\varepsilon)}\left(\sum_{j=k-w}^{L-1}\|(Tf_j)\chi_k\|_{L^{q(\cdot)}(\mathbb{R}^{n})}\right)^{p(1+\varepsilon)}\right]^{\frac{1}{p(1+\varepsilon)}}\\
&\displaystyle\hspace{3mm}+\sup_{\varepsilon>0}\sup_{L\leq0,L\in\mathbb{Z}}b^{-L\lambda}\left[\varepsilon^\theta\sum_{k=-\infty}^L b^{k\alpha(0)p(1+\varepsilon)}\left(\sum_{j=L}^{\infty}\|(Tf_j)\chi_k\|_{L^{q(\cdot)}(\mathbb{R}^{n})}\right)^{p(1+\varepsilon)}\right]^{\frac{1}{p(1+\varepsilon)}}\\
&\displaystyle=J_{11}+J_{12}+J_{13}.

\end{array}\eqno(4.1)$$

Let us first estimate $J_{11}$. By (3.3), (3.5) and the generalized H\"{o}lder
inequality, we get
$$\begin{array}{rl}
\displaystyle |Tf_j(x)|&\displaystyle\lesssim\rho(x)^{-1}\int_{B_j}|f_j(y)|dy\\
&\displaystyle\lesssim b^{-k}\|f_j\|_{L^{q(\cdot)}(\mathbb{R}^{n})}\|\chi_{B_j}\|_{L^{q'(\cdot)}(\mathbb{R}^{n})}.

\end{array}$$
So by Lemma 2.2 and Lemma 2.3, we have
$$\begin{array}{rl}
\displaystyle \|(Tf_j)\chi_k\|_{L^{q(\cdot)}(\mathbb{R}^{n})}&\displaystyle\lesssim b^{-k}\|f_j\|_{L^{q(\cdot)}(\mathbb{R}^{n})}\|\chi_{B_j}\|_{L^{q'(\cdot)}(\mathbb{R}^{n})}\|\chi_{B_k}\|_{L^{q(\cdot)}(\mathbb{R}^{n})}\\
&\displaystyle\lesssim b^{-k}\|f_j\|_{L^{q(\cdot)}(\mathbb{R}^{n})}\big(|B_k|\|\chi_{B_k}\|^{-1}_{L^{q'(\cdot)}(\mathbb{R}^{n})}\big)\|\chi_{B_j}\|_{L^{q'(\cdot)}(\mathbb{R}^{n})}\\
&\displaystyle\lesssim\|f_j\|_{L^{q(\cdot)}(\mathbb{R}^{n})}\frac{\|\chi_{B_j}\|_{L^{q'(\cdot)}(\mathbb{R}^{n})}}{\|\chi_{B_k}\|_{L^{q'(\cdot)}(\mathbb{R}^{n})}}\\
&\displaystyle\lesssim b^{\delta_2(j-k)}\|f_j\|_{L^{q(\cdot)}(\mathbb{R}^{n})}.

\end{array}\eqno(4.2)$$
Take $\frac{1}{p(1+\varepsilon)}+\frac{1}{[p(1+\varepsilon)]'}=1$. Since $0<\alpha(0)<\delta_2$, by
(4.2) and the H\"{o}lder inequality, we have
$$\begin{array}{rl}
\displaystyle J_{11}&\displaystyle\lesssim\sup_{\varepsilon>0}\sup_{L\leq0,L\in\mathbb{Z}}b^{-L\lambda}\left[\varepsilon^\theta\sum_{k=-\infty}^L b^{k\alpha(0)p(1+\varepsilon)}\left(\sum_{j=-\infty}^{k-w-1}b^{\delta_2(j-k)}\|f_j\|_{L^{q(\cdot)}(\mathbb{R}^{n})}\right)^{p(1+\varepsilon)}\right]^{\frac{1}{p(1+\varepsilon)}}\\
&\displaystyle\lesssim\sup_{\varepsilon>0}\sup_{L\leq0,L\in\mathbb{Z}}b^{-L\lambda}\bigg[\varepsilon^\theta\sum_{k=-\infty}^L \bigg(\sum_{j=-\infty}^{k-w-1}b^{j\alpha(0)p(1+\varepsilon)}b^{(j-k)[\delta_2-\alpha(0)]p(1+\varepsilon)/2}\|f_j\|^{p(1+\varepsilon)}_{L^{q(\cdot)}(\mathbb{R}^{n})}\bigg)\\
&\displaystyle\hspace{3mm}\times\bigg(\sum_{j=-\infty}^{k-w-1}b^{(j-k)[\delta_2-\alpha(0)][p(1+\varepsilon)]'/2}\bigg)^{\frac{p(1+\varepsilon)}{[p(1+\varepsilon)]'}}\bigg]^{\frac{1}{p(1+\varepsilon)}}\\
&\displaystyle\lesssim\sup_{\varepsilon>0}\sup_{L\leq0,L\in\mathbb{Z}}b^{-L\lambda}\left[\varepsilon^\theta\sum_{j=-\infty}^{L-w-1}b^{j\alpha(0)p(1+\varepsilon)}\sum_{k=j+w+1}^L b^{(j-k)(\delta_2-\alpha(0))p(1+\varepsilon)/2}\|f_j\|^{p(1+\varepsilon)}_{L^{q(\cdot)}(\mathbb{R}^{n})}\right]^{\frac{1}{p(1+\varepsilon)}}\\
&\displaystyle\lesssim\sup_{\varepsilon>0}\sup_{L\leq0,L\in\mathbb{Z}}b^{-L\lambda}\left[\varepsilon^\theta\sum_{j=-\infty}^{L-w-1}b^{j\alpha(0)p(1+\varepsilon)}\|f_j\|^{p(1+\varepsilon)}_{L^{q(\cdot)}(\mathbb{R}^{n})}\right]^{\frac{1}{p(1+\varepsilon)}}
\lesssim\|f\|_{M\dot{K}^{\alpha(\cdot),\lambda,\theta}_{q(\cdot),p)}(A;\mathbb{R}^{n})}.

\end{array}\eqno(4.3)$$

Let us now estimate $J_{12}$. By $\alpha(0)>0$,
$L^{q(\cdot)}(\mathbb{R}^{n})$ boundedness of $T$ and the H\"{o}lder
inequality, we have
$$\begin{array}{rl}
\displaystyle J_{12}&\displaystyle\lesssim\sup_{\varepsilon>0}\sup_{L\leq0,L\in\mathbb{Z}}b^{-L\lambda}\left[\varepsilon^\theta\sum_{k=-\infty}^L b^{k\alpha(0)p(1+\varepsilon)}\left(\sum_{j=k-w}^{L-1}\|f_j\|_{L^{q(\cdot)}(\mathbb{R}^{n})}\right)^{p(1+\varepsilon)}\right]^{\frac{1}{p(1+\varepsilon)}}\\
&\displaystyle\lesssim\sup_{\varepsilon>0}\sup_{L\leq0,L\in\mathbb{Z}}b^{-L\lambda}\bigg[\varepsilon^\theta\sum_{k=-\infty}^L\bigg(\sum_{j=k-w}^{L-1}b^{j\alpha(0)p(1+\varepsilon)} \|f_j\|^{p(1+\varepsilon)}_{L^{q(\cdot)}(\mathbb{R}^{n})} b^{(k-j)\alpha(0)p(1+\varepsilon)/2}\bigg)\\
&\displaystyle\hspace{3mm}\times\bigg(\sum_{j=k-w}^{L-1}b^{(k-j)\alpha(0)[p(1+\varepsilon)]'/2}\bigg)^{\frac{p(1+\varepsilon)}{[p(1+\varepsilon)]'}}\bigg]^{\frac{1}{p(1+\varepsilon)}}\\
&\displaystyle\lesssim\sup_{\varepsilon>0}\sup_{L\leq0,L\in\mathbb{Z}}b^{-L\lambda}\left[\varepsilon^\theta\sum_{j=-\infty}^{L-1}b^{j\alpha(0)p(1+\varepsilon)} \|f_j\|^{p(1+\varepsilon)}_{L^{q(\cdot)}(\mathbb{R}^{n})}\sum_{k=-\infty}^{j+w} b^{(k-j)\alpha(0)p(1+\varepsilon)/2}\right]^{\frac{1}{p(1+\varepsilon)}}\\\\
&\displaystyle\lesssim\|f\|_{M\dot{K}^{\alpha(\cdot),\lambda,\theta}_{q(\cdot),p)}(A;\mathbb{R}^{n})}.

\end{array}\eqno(4.4)$$

For $J_{13}$, since $0<2\lambda<\alpha(0)$, by
$L^{q(\cdot)}(\mathbb{R}^{n})$ boundedness of $T$ and the H\"{o}lder inequality, we have
$$\begin{array}{rl}
\displaystyle J_{13}&\displaystyle\lesssim\sup_{\varepsilon>0}\sup_{L\leq0,L\in\mathbb{Z}}b^{-L\lambda}\left[\varepsilon^\theta\sum_{k=-\infty}^L b^{k\alpha(0)p(1+\varepsilon)}\left(\sum_{j=L}^{\infty}\|f_j\|_{L^{q(\cdot)}(\mathbb{R}^{n})}\right)^{p(1+\varepsilon)}\right]^{\frac{1}{p(1+\varepsilon)}}\\
&\displaystyle\lesssim\sup_{\varepsilon>0}\sup_{L\leq0,L\in\mathbb{Z}}b^{-L\lambda}\bigg[\varepsilon^\theta\sum_{k=-\infty}^L\bigg(\sum_{j=L}^\infty b^{j\alpha(0)p(1+\varepsilon)} \|f_j\|^{p(1+\varepsilon)}_{L^{q(\cdot)}(\mathbb{R}^{n})} b^{(k-j)\alpha(0)p(1+\varepsilon)/2}\bigg)\\
&\displaystyle\hspace{3mm}\times\bigg(\sum_{j=L}^\infty b^{(k-j)\alpha(0)[p(1+\varepsilon)]'/2}\bigg)^{\frac{p(1+\varepsilon)}{[p(1+\varepsilon)]'}}\bigg]^{\frac{1}{p(1+\varepsilon)}}\\

\end{array}$$
$$\begin{array}{rl}
&\displaystyle\lesssim\sup_{\varepsilon>0}\sup_{L\leq0,L\in\mathbb{Z}}b^{-L\lambda}\left[\varepsilon^\theta\sum_{k=-\infty}^L \sum_{j=L}^{\infty} b^{(k-j)\alpha(0)p(1+\varepsilon)/2}b^{j\alpha(0)p(1+\varepsilon)}\|f_j\|^{p(1+\varepsilon)}_{L^{q(\cdot)}(\mathbb{R}^{n})}\right]^{\frac{1}{p(1+\varepsilon)}}\\
&\displaystyle\lesssim\sup_{\varepsilon>0}\sup_{L\leq0,L\in\mathbb{Z}}b^{-L\lambda}\left(\sum_{k=-\infty}^L \sum_{j=L}^{\infty} b^{(k-j)\alpha(0)p(1+\varepsilon)/2}b^{j\lambda p(1+\varepsilon)}\right)^{\frac{1}{p(1+\varepsilon)}}\|f\|_{M\dot{K}^{\alpha(\cdot),\lambda,\theta}_{q(\cdot),p)}(A;\mathbb{R}^{n})}\\
&\displaystyle\lesssim\sup_{L\leq0,L\in\mathbb{Z}}b^{-L\lambda}\left(\sum_{k=-\infty}^L b^{k\alpha(0)/2}\right)\left(\sum_{j=L}^{\infty} b^{(\lambda-\frac{\alpha(0)}{2}) j}\right)\|f\|_{M\dot{K}^{\alpha(\cdot),\lambda,\theta}_{q(\cdot),p)}(A;\mathbb{R}^{n})}\\
&\displaystyle\lesssim\sup_{L\leq0,L\in\mathbb{Z}}b^{-L\lambda} b^{L\alpha(0)/2} b^{(\lambda-\frac{\alpha(0)}{2}) L}\|f\|_{M\dot{K}^{\alpha(\cdot),\lambda,\theta}_{q(\cdot),p)}(A;\mathbb{R}^{n})}\\
&\displaystyle\lesssim\|f\|_{M\dot{K}^{\alpha(\cdot),\lambda,\theta}_{q(\cdot),p)}(A;\mathbb{R}^{n})}.

\end{array}\eqno(4.5)$$

Combining (4.1), (4.3)-(4.5), we have
$$J_1\lesssim\|f\|_{M\dot{K}^{\alpha(\cdot),\lambda,\theta}_{q(\cdot),p)}(A;\mathbb{R}^{n})}.$$

Thus, the proof of Theorem 4.3 is completed.

\section{Anisotropic grand Herz-type Hardy
spaces with variable exponents\label{s5}}$\indent$

In this section, we first give the definition of anisotropic grand Herz-type Hardy
spaces with two variable exponents $H\dot{K}^{\alpha(\cdot),p)}_{q(\cdot),\theta}(A;\mathbb{R}^{n})$ and
$HK^{\alpha(\cdot),p)}_{q(\cdot),\theta}(A;\mathbb{R}^{n})$.

We say that a $\mathcal{C}^\infty$ complex-valued function $\varphi$ on $\mathbb{R}^{n}$ belongs to the Schwartz class $\mathcal{S}$, if for every multi-index $\alpha$ and integer $m\geq0$ we have
$$\|\varphi\|_{\alpha,m}:=\sup_{x\in\mathbb{R}^{n}}\rho(x)^m|\partial^\alpha\varphi(x)|<\infty.$$

The dual space of $\mathcal{S}$ is denoted by $\mathcal{S'}$.
For integer $N\geq 0$, denote
$$\mathcal{S}_N=\{\varphi\in\mathcal{S}:\|\varphi\|_{\alpha,m}\leq 1 \,\,\mathrm{for} \,\,|\alpha|\leq N, m\leq N\}.$$
For $\varphi\in\mathcal{S}$ and $k\in\mathbb{Z}$, define the dilate of $\varphi$ to the scale $k$ by
$$\varphi_k(x)=b^{-k}\varphi(A^{-k}x).$$
In particular, if we take $A=2\mathrm{Id}$, where $\mathrm{Id}$ denotes the unit matrix, then the dilations associated
with $A$ are the usual isotropic dyadic dilations.

Suppose $f\in\mathcal{S'}$. The nontangential maximal function of $f$ with respect to $\varphi$ is defined as
$$M_\varphi(f)(x):=\sup\{|f\ast\varphi_k(y)|: x-y\in B_k,\,k\in\mathbb{Z}\}.$$
The radial maximal function of $f$ with respect to $\varphi$ is defined as
$$M^0_\varphi(f)(x):=\sup_{k\in\mathbb{Z}}|f\ast\varphi_k(y)|.$$
For given $N\in\mathbb{N}$, we define the nontangential grand maximal function of $f$ as
$$M_N(f)(x):=\sup_{\varphi\in\mathcal{S}_N}M_\varphi(f)(x).$$
The radial grand maximal function of $f$ is
$$M^0_N(f)(x):=\sup_{\varphi\in\mathcal{S}_N}M^0_\varphi(f)(x).$$

From [2, Proposition 3.10], we know that radial and nontangential grand maximal functions are pointwise equivalent, i.e., for every $N\geq0$, there is a constant $C=C(N)$ so that, for all $f\in\mathcal{S'}$ and $x\in \mathbb{R}^{n}$,
$$M^0_N(f)(x)\leq M_N(f)(x)\leq CM^0_N(f)(x).\eqno(5.1)$$

Now let us introduce the anisotropic grand Herz-type Hardy spaces with two variable exponents as follows.

\noindent{\bf Definition 5.1}\quad Let $\alpha(\cdot):
\mathbb{R}^{n}\rightarrow\mathbb{R}$ with $\alpha(\cdot)\in L^\infty(\mathbb{R}^{n})$, $1\leq p<\infty$, $q(\cdot)\in
\mathcal{P}(\mathbb{R}^{n})$, $\theta>0$ and $N>2$.

(i) The homogeneous anisotropic grand Herz-type Hardy space
$H\dot{K}^{\alpha(\cdot),p)}_{q(\cdot),\theta}(A;\mathbb{R}^{n})$ associated with the dilation $A$ is defined by
$$H\dot{K}^{\alpha(\cdot),p)}_{q(\cdot),\theta}(A;\mathbb{R}^{n})=\{f\in \mathcal{S'}(\mathbb{R}^{n}):
M_Nf(x)\in\dot{K}^{\alpha(\cdot),p)}_{q(\cdot),\theta}(A;\mathbb{R}^{n})\}$$ and we
define
$\|f\|_{H\dot{K}^{\alpha(\cdot),p)}_{q(\cdot),\theta}(A;\mathbb{R}^{n})}=\|M_Nf\|_{\dot{K}^{\alpha(\cdot),p)}_{q(\cdot),\theta}(A;\mathbb{R}^{n})}$.

(ii) The non-homogeneous anisotropic grand Herz-type Hardy space $HK^{\alpha(\cdot),p)}_{q(\cdot),\theta}(A;\mathbb{R}^{n})$ associated with the dilation $A$ is defined by
$$HK^{\alpha(\cdot),p)}_{q(\cdot),\theta}(A;\mathbb{R}^{n})=\{f\in \mathcal{S'}(\mathbb{R}^{n}):
M_Nf(x)\in K^{\alpha(\cdot),p)}_{q(\cdot),\theta}(A;\mathbb{R}^{n})\}$$ and we define
$\|f\|_{HK^{\alpha(\cdot),p)}_{q(\cdot),\theta}(A;\mathbb{R}^{n})}=\|M_Nf\|_{K^{\alpha(\cdot),p)}_{q(\cdot),\theta}(A;\mathbb{R}^{n})}$.

It is obvious that $M_Nf$ satisfies (3.11). Thus, by Theorem 3.3, we
can easily prove that if $\alpha(\cdot)\in L^\infty(\mathbb{R}^{n})\cap\mathcal{P}_0(\mathbb{R}^{n})\cap\mathcal{P}_\infty(\mathbb{R}^{n})$, $0<\alpha(0), \alpha_\infty<\delta_2$, $1\leq p<\infty$ and
$q(\cdot)\in \mathcal{B}(\mathbb{R}^{n})$, then
$$H\dot{K}^{\alpha(\cdot),p)}_{q(\cdot),\theta}(A;\mathbb{R}^{n})\cap L^{q(\cdot)}_{\mathrm{loc}}(\mathbb{R}^{n}\setminus\{0\})=\dot{K}^{\alpha(\cdot),p)}_
{q(\cdot),\theta}(A;\mathbb{R}^{n})$$
and $$HK^{\alpha(\cdot),p)}_{q(\cdot),\theta}(A;\mathbb{R}^{n})\cap
L^{q(\cdot)}_{\mathrm{loc}}(\mathbb{R}^{n})=K^{\alpha(\cdot),p)}_{q(\cdot),\theta}(A;\mathbb{R}^{n}).$$
If $\alpha(\cdot)\in L^\infty(\mathbb{R}^{n})\cap\mathcal{P}_0(\mathbb{R}^{n})\cap\mathcal{P}_\infty(\mathbb{R}^{n})$, $\delta_2\leq\alpha(0), \alpha_\infty<\infty$, $1\leq p<\infty$ and $q(\cdot)\in
\mathcal{B}(\mathbb{R}^{n})$, then
$$H\dot{K}^{\alpha(\cdot),p)}_{q(\cdot),\theta}(A;\mathbb{R}^{n})\cap L^{q(\cdot)}_{\mathrm{loc}}(\mathbb{R}^{n}\setminus\{0\})\subsetneqq\dot{K}^
{\alpha(\cdot),p)}_{q(\cdot),\theta}(A;\mathbb{R}^{n})$$
and $$HK^{\alpha(\cdot),p)}_{q(\cdot),\theta}(A;\mathbb{R}^{n})\cap
L^{q(\cdot)}_{\mathrm{loc}}(\mathbb{R}^{n})\subsetneqq
K^{\alpha(\cdot),p)}_{q(\cdot),\theta}(A;\mathbb{R}^{n}).$$ Thus we are interested in
the case $\alpha\geq \delta_2$. In this case, we establish
characterizations of the spaces
$H\dot{K}^{\alpha(\cdot),p)}_{q(\cdot),\theta}(A;\mathbb{R}^{n})$ and
$HK^{\alpha(\cdot),p)}_{q(\cdot),\theta}(A;\mathbb{R}^{n})$ in terms of central
atomic decompositions. For $x\in \mathbb{R}$ we denote by $[x]$ the
largest integer less than or equal to $x$.

\noindent{\bf Definition 5.2}\quad Let $\alpha(\cdot)\in L^\infty(\mathbb{R}^{n})\cap\mathcal{P}_0(\mathbb{R}^{n})\cap\mathcal{P}_\infty(\mathbb{R}^{n})$, $\delta_2\leq\alpha(0), \alpha_\infty<\infty$,
$q(\cdot)\in \mathcal{P}(\mathbb{R}^{n})$, and non-negative integer
$s\geq [(\alpha-\delta_2)\ln b/\ln\lambda_{-}]$.

(i) A function $a$ on $\mathbb{R}^{n}$ is said to be a central
$(\alpha(\cdot), q(\cdot))$-atom, if it satisfies

\hspace{3mm}(1) supp\,$a\subset B_k$, $k\in\mathbb{Z}$.

\hspace{3mm}(2) $\|a\|_{L^{q(\cdot)}(\mathbb{R}^{n})}\leq
|B_k|^{-\alpha(0)}$, $k<0$.

\hspace{3mm}(3) $\|a\|_{L^{q(\cdot)}(\mathbb{R}^{n})}\leq
|B_k|^{-\alpha_\infty}$, $k\geq 0$.

\hspace{3mm}(4) $\int_{\mathbb{R}^{n}}a(x)x^\beta dx=0, |\beta|\leq
s$.

(ii) A function $a$ on $\mathbb{R}^{n}$ is said to be a central
$(\alpha(\cdot), q(\cdot))$-atom of restricted type, if it satisfies the
conditions (3), (4) above and

\hspace{3mm}(1)$'$ supp\,$a\subset B_k, k\geq 0$.

\noindent {\bf Remark 5.1}\quad If $\alpha(x)=\alpha$ is a constant, then we can get the atomic definition in [39]. If $\alpha(x)=\alpha$ and $q(x)=q$ are constants, then taking $\delta_2=1-1/q$ we can get the classical case in [9].

Next we give the atomic decomposition theorems.

\noindent{\bf Theorem 5.1}\quad Let $\alpha(\cdot)\in L^\infty(\mathbb{R}^{n})\cap\mathcal{P}_0(\mathbb{R}^{n})\cap\mathcal{P}_\infty(\mathbb{R}^{n})$, $\delta_2\leq\alpha(0), \alpha_\infty<\delta_2+\ln\lambda_{-}/\ln b$,
$1\leq p<\infty$ and $q(\cdot)\in \mathcal{B}(\mathbb{R}^{n})$. Then we
have

(i) $f\in H\dot{K}^{\alpha(\cdot),p)}_{q(\cdot),\theta}(A;\mathbb{R}^{n})$ if and
only if
$$f=\sum_{k=-\infty}^\infty\lambda_ka_k,\quad \mathrm{in\,\, the\,\, sense\,\, of\,\,} \mathcal{S'}(\mathbb{R}^{n}),\eqno(5.2)$$ where each
$a_k$ is a central $(\alpha(\cdot), q(\cdot))$-atom with support contained
in $B_k$ and
$\displaystyle\sup_{\varepsilon>0}\varepsilon^\theta\sum_{k=-\infty}^\infty|\lambda_k|^{p(1+\varepsilon)}<\infty$.
Moreover,
$$\|f\|_{H\dot{K}^{\alpha(\cdot),p)}_{q(\cdot),\theta}(A;\mathbb{R}^{n})}\approx\inf\displaystyle\sup_{\varepsilon>0}\left(\varepsilon^\theta\sum_{k=-\infty}^\infty|\lambda_k|^{p(1+\varepsilon)}\right)^{\frac{1}{p(1+\varepsilon)}}<\infty,$$
where the infimum is taken over all above decompositions of $f$.

(ii) $f\in HK^{\alpha(\cdot),p)}_{q(\cdot),\theta}(A;\mathbb{R}^{n})$ if and only if
$$f=\sum_{k=0}^\infty\lambda_ka_k,\quad \rm{in\,\, the\,\, sense\,\, of\,\,} \mathcal{S'}(\mathbb{R}^{n}),$$ where each
$a_k$ is a central $(\alpha(\cdot), q(\cdot))$-atom of restricted type with
support contained in $B_k$ and
$\displaystyle\sup_{\varepsilon>0}\varepsilon^\theta\sum_{k=0}^\infty|\lambda_k|^{p(1+\varepsilon)}<\infty$. Moreover,
$$\|f\|_{HK^{\alpha(\cdot),p)}_{q(\cdot),\theta}(A;\mathbb{R}^{n})}\approx\inf\displaystyle\sup_{\varepsilon>0}\left(\varepsilon^\theta\sum_{k=0}^\infty|\lambda_k|^{p(1+\varepsilon)}\right)^{\frac{1}{p(1+\varepsilon)}}<\infty,$$
where the infimum is taken over all above decompositions of $f$.

\noindent{\bf Proof}\quad We only prove (i), and (ii) can be proved
in the similar way. To prove the necessity, choose $\varphi\in
\mathcal{S}(\mathbb{R}^n)$ such that $\int_{\mathbb{R}^n}\varphi(x)dx=1$. Set
$f^{(j)}(x):=f\ast \varphi_{-j}(x)$ for each $f\in H\dot{K}^{\alpha(\cdot),p)}_{q(\cdot),\theta}(A;\mathbb{R}^{n})$ and $j\in \mathbb{N}$.
It is obvious that $f^{(j)}\in C^\infty(\mathbb{R}^{n})$ and
$\displaystyle\lim_{j\rightarrow\infty}f^{(j)}=f$ in the sense of
distribution. Let $\psi\in C^\infty_0(\mathbb{R}^n)$ such that
supp\,$\psi\subset \tilde{C}_0:=C_{-1}\cup C_0\cup C_1$, $0\leq \psi\leq 1$, and
$\psi(x)=1$ for $x\in C_0$. Let
$\psi_{(k)}(x)=\psi(A^{-k}x)$ for $k\in \mathbb{Z}$.
Obviously, supp\,$\psi_{(k)}\subset \tilde{C}_k=C_{k-1}\cup C_k\cup C_{k+1}$. Let$$\displaystyle \Psi_k(x)=\left\{\begin{array}{ll}
\displaystyle \psi_{(k)}(x)\bigg/\left(\sum_{l=-\infty}^\infty\psi_l(x)\right),\quad x\neq 0,\\
\displaystyle 0,\hspace{4cm} x=0,
\end{array}\right.$$
then we have $\Psi_k\in C^\infty_0(\mathbb{R}^n)$, supp\,$\Psi_{k}\subset \tilde{C}_k$, $0\leq \Psi_k\leq 1$, and $\displaystyle\sum_{k=-\infty}^\infty\Psi_k(x)=1$ for $x\neq
0$. Set $\nu_k(x)=|\tilde{C}_k|^{-1}\chi_{\tilde{C}_{k}}(x)$. We write
$$\begin{array}{rl}
\displaystyle f^{(j)}(x)
&\displaystyle=\sum_{k=-\infty}^\infty f^{(j)}(x)\Psi_k(x)\\
&\displaystyle=\sum_{k=-\infty}^\infty \left\{f^{(j)}(x)\Psi_k(x)-\left(\int_{\mathbb{R}^{n}}f^{(j)}(y)\Psi_k(y)dy\right)\nu_k(x)\right\}\\
&\displaystyle\hspace{3mm}+\sum_{k=-\infty}^\infty\left(\int_{\mathbb{R}^{n}}f^{(j)}(y)\Psi_k(y)dy\right)\nu_k(x)\\
&\displaystyle:=\sum_I^{(j)}+\sum_{II}^{(j)}.

\end{array}$$

For the term $\displaystyle\sum_I^{(j)}$, let
$$g_k^{(j)}(x)=f^{(j)}(x)\Psi_k(x)-\left(\int_{\mathbb{R}^{n}}f^{(j)}(y)\Psi_k(y)dy\right)\nu_k(x)$$ and
$a_{1,k}^{(j)}(x)=g_k^{(j)}(x)/\lambda_{1,k}$, where
$\displaystyle\lambda_{1,k}=b_1\sum_{l=k-1}^{k+1}\|b^{(k+1)\alpha(\cdot)}(M_Nf)\chi_l\|_{L^{q(\cdot)}(\mathbb{R}^{n})}$,
and $b_1$ is a constant which will be chosen later. Note that
supp\,$\displaystyle a_{1,k}^{(j)}\subset \tilde{C}_k\subset B_{k+1},\,\int_{\mathbb{R}^{n}}a_{1,k}^{(j)}(x)dx=0$, and
$\displaystyle\sum_I^{(j)}=\sum_{k=-\infty}^\infty\lambda_{1,k}a_{1,k}^{(j)}(x)$.

Now we estimate $\|b^{(k+1)\alpha(\cdot)}g_k^{(j)}\|_{L^{q(\cdot)}(\mathbb{R}^{n})}$. By the generalized H\"{o}lder
inequality and Lemma 2.2 we have
$$\begin{array}{rl}
&\displaystyle\|b^{(k+1)\alpha(\cdot)}g_k^{(j)}\|_{L^{q(\cdot)}(\mathbb{R}^{n})}\\
&\displaystyle\leq
\|b^{(k+1)\alpha(\cdot)}f^{(j)}\Psi_k\|_{L^{q(\cdot)}(\mathbb{R}^{n})}+\frac{C}{|\tilde{C}_k|}\|b^{(k+1)\alpha(\cdot)}f^{(j)}\Psi_k\|_{L^{q(\cdot)}(\mathbb{R}^{n})}
\|\chi_{\tilde{C}_k}\|_{L^{q'(\cdot)}(\mathbb{R}^{n})}\|\chi_{\tilde{C}_k}\|_{L^{q(\cdot)}(\mathbb{R}^{n})}\\
&\displaystyle\leq
\|b^{(k+1)\alpha(\cdot)}f^{(j)}\Psi_k\|_{L^{q(\cdot)}(\mathbb{R}^{n})}+C\|b^{(k+1)\alpha(\cdot)}f^{(j)}\Psi_k\|_{L^{q(\cdot)}(\mathbb{R}^{n})}\\
&\displaystyle\leq
C'\sum_{l=k-1}^{k+1}\|b^{(k+1)\alpha(\cdot)}(M_Nf)\chi_l\|_{L^{q(\cdot)}(\mathbb{R}^{n})}.

\end{array}$$
Choose $b_1=C'$, then
$\|b^{(k+1)\alpha(\cdot)}a_{1,k}^{(j)}\|_{L^{q(\cdot)}(\mathbb{R}^{n})}\leq
C$ and each $a_{1,k}^{(j)}$ is a central $(\alpha(\cdot),
q(\cdot))$-atom with support contained in $B_{k+1}$ by Definition 5.2. Furthermore,
$$\begin{array}{rl}
\displaystyle\sup_{\varepsilon>0}\varepsilon^\theta\sum_{k=-\infty}^\infty|\lambda_{1,k}|^{p(1+\varepsilon)}&\displaystyle\leq C\sup_{\varepsilon>0}\varepsilon^\theta\sum_{k=-\infty}^\infty\left(\sum_{l=k-1}^{k+1}\|b^{(k+1)\alpha(\cdot)}(M_Nf)\chi_l\|_{L^{q(\cdot)}(\mathbb{R}^{n})}\right)^{p(1+\varepsilon)}\\
&\displaystyle\leq C\sup_{\varepsilon>0}\varepsilon^\theta\|M_Nf\|^{p(1+\varepsilon)}_{\dot{K}^{\alpha(\cdot),p(1+\varepsilon)}_{q(\cdot)}(A;\mathbb{R}^{n})}\\
&\displaystyle\leq C\|M_Nf\|^{p(1+\varepsilon)}_{\dot{K}^{\alpha(\cdot),p),\theta}_{q(\cdot)}(A;\mathbb{R}^{n})}<\infty.

\end{array}$$

It remains to estimate $\displaystyle\sum_{II}^{(j)}$. Summing by parts, we have
$$\displaystyle\sum_{II}^{(j)}=\sum_{k=-\infty}^\infty\left(\sum_{j=-\infty}^k\int_{\mathbb{R}^{n}}f^{(j)}(y)\Psi_j(y)dy\right)(\nu_k(x)-\nu_{k+1}(x)):=\sum_{k=-\infty}^\infty h^{(j)}_k(x),$$
where $$h^{(j)}_k(x)=\left(\sum_{j=-\infty}^k\int_{\mathbb{R}^{n}}f^{(j)}(y)\Psi_j(y)dy\right)(\nu_k(x)-\nu_{k+1}(x)).$$
Let
$a_{2,k}^{(j)}(x)=h_k^{(j)}(x)/\lambda_{2,k}$, where
$\displaystyle\lambda_{2,k}=b_2\sum_{l=k-1}^{k+2}\|b^{(k+2)\alpha(\cdot)}(M_Nf)\chi_l\|_{L^{q(\cdot)}(\mathbb{R}^{n})}$,
and $b_2$ is a constant which will be chosen later. Note that
supp\,$\displaystyle a_{2,k}^{(j)}\subset \tilde{C}_{k+1}\subset B_{k+2},\,\int_{\mathbb{R}^{n}}a_{2,k}^{(j)}(x)dx=0$, and
$\displaystyle\sum_{II}^{(j)}=\sum_{k=-\infty}^\infty\lambda_{2,k}a_{2,k}^{(j)}(x)$.
Let $\displaystyle \Phi(x)=\sum_{l=-\infty}^{-2}\Psi_l(x)$. Since supp\,$\displaystyle \Psi_l\subset \tilde{C}_{l}$ and $\{\tilde{C}_{l}\}_{}l=-\infty^{-2}$ has bounded overlapping, that is $\displaystyle \sum_{l=-\infty}^{-2}\chi_{\tilde{C}_{l}}\leq C$, it is easy to see that $\Phi\in C_0^\infty(B_{-1})$ and hence $\Phi\in\mathcal{S}$. Note that
$$\sum_{l=-\infty}^{k}\Psi_l(x)=\Phi(A^{-k-2}x)=b^{k+2}\Phi_{k+2}(x),$$
thus for any $x\in B_{k+2}$, by using [4, Lemma 6.6], we obtain
$$\begin{array}{rl}
\displaystyle\left|\int_{\mathbb{R}^{n}}f^{(j)}(y)\sum_{j=-\infty}^k\Psi_j(y)dy\right|
&\displaystyle=b^{k+2}\left|\int_{B^{k+2}}f^{(j)}(y)\Phi_{k+2}(y)dy\right|\\
&\displaystyle=b^{k+2}\left|f^{(j)}\ast\tilde{\Phi}_{k+2}(0)\right|\\
&\displaystyle\leq b^{k+2}\|\tilde{\Phi}\|_{\mathcal{S}_{N+2}}M_{N+2}f^{(j)}(x)\\
&\displaystyle\leq Cb^{k+2}M_{N}f(x),

\end{array}$$
where
$\tilde{\Phi}(y)=\Phi(-y)$ and $C$ is a constant dependent of $N$. On the other hand,
$$|\nu_k(x)-\nu_{k+1}(x)|\leq Cb^{-k-2}\sum_{l=k-1}^{k+2}\chi_l(x).$$
Hence,
$$\|b^{(k+2)\alpha(\cdot)}h_k^{(j)}\|_{L^{q(\cdot)}(\mathbb{R}^{n})}\leq C''\sum_{l=k-1}^{k+2}\|b^{(k+2)\alpha(\cdot)}(M_Nf)\chi_l\|_{L^{q(\cdot)}(\mathbb{R}^{n})}.$$
Take $b_2=C''$. It is easy to verify that each
$a_{2,k}^{(j)}$ is a central $(\alpha(\cdot), q(\cdot))$-atom with support
contained in
$\tilde{C}_{k}\cup\tilde{C}_{k+1}\subset
B_{k+2}$.
Moreover,
$$\begin{array}{rl}
\displaystyle\sup_{\varepsilon>0}\varepsilon^\theta\sum_{k=-\infty}^\infty|\lambda_{2,k}|^{p(1+\varepsilon)}&\displaystyle\leq C\sup_{\varepsilon>0}\varepsilon^\theta\sum_{k=-\infty}^\infty\left(\sum_{l=k-1}^{k+1}\|b^{(k+2)\alpha(\cdot)}(M_Nf)\chi_l\|_{L^{q(\cdot)}(\mathbb{R}^{n})}\right)^{p(1+\varepsilon)}\\
&\displaystyle\leq C\sup_{\varepsilon>0}\varepsilon^\theta\|M_Nf\|^{p(1+\varepsilon)}_{\dot{K}^{\alpha(\cdot),p(1+\varepsilon)}_{q(\cdot)}(A;\mathbb{R}^{n})}\\
&\displaystyle\leq C\|M_Nf\|^{p(1+\varepsilon)}_{\dot{K}^{\alpha(\cdot),p),\theta}_{q(\cdot)}(A;\mathbb{R}^{n})}<\infty.

\end{array}$$

Thus we obtain that
$$f^{(j)}(x)=\sum_{d=-\infty}^\infty \lambda_da_d^{(j)}(x),$$ where
each $a_d^{(j)}$ is a central $(\alpha(\cdot),q(\cdot))$-atom with support
contained in
$B_{d+2}$, $\lambda_d$ is independent of $j$ and
$$\sup_{\varepsilon>0}\left(\varepsilon^\theta\sum_{d=-\infty}^\infty|\lambda_d|^{p(1+\varepsilon)}\right)^{\frac{1}{p(1+\varepsilon)}}\leq C \|M_Nf\|_{\dot{K}^{\alpha(\cdot),p),\theta}_{q(\cdot)}(A;\mathbb{R}^{n})}<\infty,$$
where $C$ is independent of $j$ and $f$.

Since $$\sup_{j\in \mathbb{Z}_+}
\|a_0^{(j)}\|_{L^{q(\cdot)}(\mathbb{R}^{n})}\leq
|B_2|^{-\alpha/n},$$ by the Banach-Alaoglu theorem we can obtain a
subsequence $\{a_0^{(j_{n_0})}\}$ of $\{a_0^{(j)}\}$ converging in
the $\mathrm{weak}^\ast$ topology of $L^{q(\cdot)}(\mathbb{R}^{n})$
to some $a_0\in L^{q(\cdot)}(\mathbb{R}^{n})$. It is readily to
verify that $a_0$ is a central $(\alpha(\cdot),q(\cdot))$-atom supported on
$B_2$. Next, since
$$\sup_{j_{n_0}\in \mathbb{Z}_+}
\|a_1^{(j_{n_0})}\|_{L^{q(\cdot)}(\mathbb{R}^{n})}\leq
|B_3|^{-\alpha},$$ another application of the Banach-Alaoglu
theorem yields a subsequence $\{a_1^{(j_{n_1})}\}$ of
$\{a_1^{(j_{n_0})}\}$ which converges $\mathrm{weak}^\ast$ in
$L^{q(\cdot)}(\mathbb{R}^{n})$ to a central $(\alpha(\cdot),q(\cdot))$-atom
$a_1$ with support in $B_3$. Furthermore, $$\sup_{j_{n_1}\in
\mathbb{Z}_+}
\|a_{-1}^{(j_{n_1})}\|_{L^{q(\cdot)}(\mathbb{R}^{n})}\leq
|B_1|^{-\alpha}.$$ Similarly, there exists a subsequence
$\{a_{-1}^{(j_{n_{-1}})}\}$ of $\{a_{-1}^{(j_{n_1})}\}$ which
converges $\mathrm{weak}^\ast$ in $L^{q(\cdot)}(\mathbb{R}^{n})$ to
some $a_{-1}\in L^{q(\cdot)}(\mathbb{R}^{n})$, and $a_{-1}$ is a
central $(\alpha(\cdot),q(\cdot))$-atom supported on $B_1$. Repeating the
above procedure for each $d\in \mathbb{Z}$, we can find a
subsequence $\{a_d^{(j_{n_d})}\}$ of $\{a_d^{(j)}\}$ converging
$\mathrm{weak}^\ast$ in $L^{q(\cdot)}(\mathbb{R}^{n})$ to some
$a_d\in L^{q(\cdot)}(\mathbb{R}^{n})$ which is a central
$(\alpha(\cdot),q(\cdot))$-atom supported on $B_{d+2}$. By the usual
diagonal method we obtain a subsequence $\{j_\nu\}$ of
$\mathbb{Z}_+$ such that for each $\displaystyle d\in
\mathbb{Z},\,\lim_{\nu\rightarrow\infty}a_d^{(j_\nu)}=a_d$ in the
$\mathrm{weak}^\ast$ topology of $L^{q(\cdot)}(\mathbb{R}^{n})$ and
therefore in $\mathcal{S'}(\mathbb{R}^{n})$.

Now our proof is reduced to prove that
$$f=\sum_{d=-\infty}^\infty\lambda_da_d, \,\,\,
\mathrm{in\,\, the\,\, sense\,\, of\,\,}
\mathcal{S'}(\mathbb{R}^{n}).\eqno(5.3)$$ For each
$\phi\in\mathcal{S}(\mathbb{R}^{n})$, note that
$\mathrm{supp}\,a_d^{(j_\nu)}\subset(\tilde{C}_{d}\cup\tilde{C}_{d+1})\subset(C_{d-1}\cup
C_d\cup C_{d+1}\cup C_{d+2})$. We have
$$\langle f, \phi\rangle=\lim_{\nu\rightarrow\infty}\sum_{d=-\infty}^\infty\lambda_d\int_{\mathbb{R}^n}
a_d^{(j_\nu)}(x)\phi(x)dx.$$ See [9] for the details.

If $d+1\leq 0$, then by Lemma 2.3
and the generalized H\"{o}lder inequality we have
$$\begin{array}{rl}
\displaystyle \left|\int_{\mathbb{R}^n} a_d^{(j_\nu)}(x)\phi(x)dx\right|&\displaystyle=\left|\int_{B_{d+2}}
a_d^{(j_\nu)}(x)\left(\phi(x)-\phi(0)\right)dx\right|\\
&\displaystyle\leq C \int_{B_{d+2}} |a_d^{(j_\nu)}(x)|\cdot |x|\sup_{z\in B_{d+2}}\sup_{|\alpha|=1}|\partial^\alpha\phi(z)|dx\\
&\displaystyle\leq C b^{(d+1)\ln\lambda_{-}/\ln b}\int_{B_{d+2}} |a_d^{(j_\nu)}(x)|dx\\
&\displaystyle\leq C
b^{(d+1)(\ln\lambda_{-}/\ln b-\alpha(0))}\|\chi_{B_{d+2}}\|_{L^{q'(\cdot)}(\mathbb{R}^{n})}\\
&\displaystyle\leq C
b^{(d+1)(\ln\lambda_{-}/\ln b-\alpha(0))}\left(\frac{|B_{d+2}|}{|B_2|}\right)^{\delta_2}\|\chi_{B_2}\|_{L^{q'(\cdot)}(\mathbb{R}^{n})}\\
&\displaystyle\leq C
b^{(d+1)(\ln\lambda_{-}/\ln b-\alpha(0)+\delta_2)}\frac{|B_2|}{|B_0|}\|\chi_{B_0}\|_{L^{q'(\cdot)}(\mathbb{R}^{n})}\\
&\displaystyle\leq C
b^{(d+1)(\ln\lambda_{-}/\ln b-\alpha(0)+\delta_2)}\inf\left\{\lambda>0: \int_{B_0}\lambda^{-q'(x)}dx\leq 1\right\}\\
&\displaystyle\leq C
b^{(d+1)(\ln\lambda_{-}/\ln b-\alpha(0)+\delta_2)}\inf\left\{1\geq\lambda>0: \int_{B_0}\lambda^{-(q')^+}dx\leq 1\right\}\\
&\displaystyle\leq C b^{d(\ln\lambda_{-}/\ln b-\alpha(0)+\delta_2)}|B_0|^{1/(q')^+}\\
&\displaystyle= C b^{d(\ln\lambda_{-}/\ln b-\alpha(0)+\delta_2)},
\end{array}$$
where $C$ is independent of $d$.

If $d+1>0$, let $k_0\in\mathbb{Z}_+$ such that $k_0+\alpha-1>0$, then
by Lemma 2.3 and the generalized H\"{o}lder inequality we have
$$\begin{array}{rl}
\displaystyle \left|\int_{\mathbb{R}^n} a_d^{(j_\nu)}(x)\phi(x)dx\right|
&\displaystyle\leq C \int_{\mathbb{R}^n}|a_d^{(j_\nu)}(x)|\rho(x)^{-k_0}dx\\
&\displaystyle\leq C
b^{-d(k_0+\alpha_\infty)}\|\chi_{B_{d+2}}\|_{L^{q'(\cdot)}(\mathbb{R}^{n})}\\
&\displaystyle\leq C
b^{-d(k_0+\alpha_\infty)}\frac{|B_{d+2}|}{|B_0|}\|\chi_{B_0}\|_{L^{q'(\cdot)}(\mathbb{R}^{n})}\\
&\displaystyle\leq C b^{-d(k_0+\alpha_\infty-1)}\inf\left\{\lambda>0: \int_{B_0}\lambda^{-q'(x)}dx\leq 1\right\}\\
&\displaystyle\leq C b^{-d(k_0+\alpha_\infty-1)}\inf\left\{1\geq\lambda>0: \int_{B_0}\lambda^{-(q')^+}dx\leq 1\right\}\\
&\displaystyle= C b^{-d(k_0+\alpha_\infty-1)},
\end{array}$$
where $C$ is independent of $d$.

Let $$\displaystyle \mu_d=\left\{\begin{array}{ll}
\displaystyle |\lambda_d|b^{d(\ln\lambda_{-}/\ln b-\alpha(0)+\delta_2)},\quad d+1\leq 0,\\
\displaystyle |\lambda_d|b^{-d(k_0+\alpha_\infty-1)},\quad\quad d+1>0,
\end{array}\right.$$
Then $$\sum_{d=-\infty}^\infty\mu_d\leq
C\sup_{\varepsilon>0}\left(\varepsilon^\theta\sum_{d=-\infty}^\infty|\lambda_d|^{p(1+\varepsilon)}\right)^{\frac{1}{p(1+\varepsilon)}}\leq C
\|M_Nf\|_{\dot{K}^{\alpha(\cdot),p),\theta}_{q(\cdot)}(A;\mathbb{R}^{n})}<\infty$$
and
$$|\lambda_d|\left|\int_{\mathbb{R}^{n}} a_d^{(j_\nu)}(x)\varphi(x)dx\right|\leq
C|\mu_d|,$$ which implies that
$$\langle f, \phi\rangle=\sum_{d=-\infty}^\infty\lim_{\nu\rightarrow\infty} \lambda_d\int_{\mathbb{R}^n} a_d^{(j_\nu)}(x)\phi(x)dx
=\sum_{d=-\infty}^\infty \lambda_d\int_{\mathbb{R}^n}
a_d(x)\phi(x)dx.$$

This establishes the identity (5.3).

To prove the sufficiency, we consider the two cases $k<0$ and
$k\geq 0$. When $k<0$, we have
$$\begin{array}{rl}
\displaystyle
\sum_{k=-\infty}^{-1}\|b^{k\alpha(\cdot)}(M^0_Nf)\chi_k\|^{p(1+\varepsilon)}_{L^{q(\cdot)}(\mathbb{R}^{n})}
&\displaystyle\leq
\sum_{k=-\infty}^{-1}b^{k\alpha(0)p(1+\varepsilon)}\left(\sum_{l=-\infty}^\infty|\lambda_l|\|(M^0_Na_l)\chi_k\|_{L^{q(\cdot)}(\mathbb{R}^{n})}\right)^{p(1+\varepsilon)}\\
&\displaystyle\leq C\sum_{k=-\infty}^{-1}b^{k\alpha(0)p(1+\varepsilon)}\left(\sum_{l=-\infty}^{k-w-1}|\lambda_l|\|(M^0_Na_l)\chi_k\|_{L^{q(\cdot)}(\mathbb{R}^{n})}\right)^{p(1+\varepsilon)}\\
&\displaystyle\hspace{3mm}+ C\sum_{k=-\infty}^{-1}b^{k\alpha(0)p(1+\varepsilon)}\left(\sum_{l=k-w}^{-1}|\lambda_l|\|(M^0_Na_l)\chi_k\|_{L^{q(\cdot)}(\mathbb{R}^{n})}\right)^{p(1+\varepsilon)}\\
&\displaystyle\hspace{3mm}+ C\sum_{k=-\infty}^{-1}b^{k\alpha(0)p(1+\varepsilon)}\left(\sum_{l=0}^\infty|\lambda_l|\|(M^0_Na_l)\chi_k\|_{L^{q(\cdot)}(\mathbb{R}^{n})}\right)^{p(1+\varepsilon)}\\
&\displaystyle=V_1+V_2+V_3.

\end{array}$$
Using the $L^{q(\cdot)}(\mathbb{R}^{n})$ boundedness of the radial grand
maximal operator $M^0_N$ and the H\"{o}lder inequality, we have
$$\begin{array}{rl}
\displaystyle V_3&\displaystyle\leq
C\sum_{k=-\infty}^{-1}b^{k\alpha(0)p(1+\varepsilon)}\left(\sum_{l=0}^\infty|\lambda_l|\|a_l\|_{L^{q(\cdot)}(\mathbb{R}^{n})}\right)^{p(1+\varepsilon)}\\
&\displaystyle\leq C\sum_{k=-\infty}^{-1}b^{k\alpha(0)p(1+\varepsilon)}\left(\sum_{l=0}^\infty|\lambda_l|^{p(1+\varepsilon)}|B_l|^{-\alpha_\infty
p(1+\varepsilon)/2}\right)\left(\sum_{l=0}^\infty|B_l|^{-\alpha_\infty
[p(1+\varepsilon)]'/2}\right)^{\frac{p(1+\varepsilon)}{[p(1+\varepsilon)]'}}\\
&\displaystyle\leq C\sum_{k=-\infty}^{-1}b^{k\alpha(0)p(1+\varepsilon)}\left(\sum_{l=0}^\infty|\lambda_l|^{p(1+\varepsilon)}\right)\leq C\sum_{l=0}^\infty|\lambda_l|^{p(1+\varepsilon)}.

\end{array}$$
For $V_2$, we have
$$\begin{array}{rl}
\displaystyle V_2&\displaystyle\leq
C\sum_{k=-\infty}^{-1}b^{k\alpha(0)p(1+\varepsilon)}\left(\sum_{l=k-w}^{-1}|\lambda_l|\|a_l\|_{L^{q(\cdot)}(\mathbb{R}^{n})}\right)^{p(1+\varepsilon)}\\
&\displaystyle\leq C\sum_{k=-\infty}^{-1}b^{k\alpha(0)p(1+\varepsilon)/2}\left(\sum_{l=k-w}^{-1}|\lambda_l|^{p(1+\varepsilon)}b^{-l\alpha(0)
p(1+\varepsilon)/2}\right)\left(\sum_{l=k-w}^{-1}b^{(k-l)\alpha(0)
[p(1+\varepsilon)]'/2}\right)^{\frac{p(1+\varepsilon)}{[p(1+\varepsilon)]'}}\\

\end{array}$$
$$\begin{array}{rl}
&\displaystyle\leq C\sum_{k=-\infty}^{-1}b^{k\alpha(0)p(1+\varepsilon)/2}\left(\sum_{l=k-w}^{-1}|\lambda_l|^{p(1+\varepsilon)}b^{-l\alpha(0)
p(1+\varepsilon)/2}\right)\\
&\displaystyle\leq C\left(\sum_{l=-\infty}^{-1}|\lambda_l|^{p(1+\varepsilon)}\right)\left(\sum_{k=-\infty}^{l+w}b^{(k-l)\alpha(0)p(1+\varepsilon)/2}\right)\\
&\displaystyle\leq C\sum_{l=-\infty}^{-1}|\lambda_l|^{p(1+\varepsilon)}.

\end{array}$$
For $V_1$, noting that $k\geq l+w+1$, similar to the method of [9], by taking $s=[(\alpha-\delta_2)\ln b/\ln\lambda_{-}]$ we have
$$\begin{array}{rl}
\displaystyle \|(M^0_Na_l)\chi_k\|_{L^{q(\cdot)}(\mathbb{R}^{n})}&\displaystyle\leq
Cb^{-l\alpha(0)-l}(b\lambda_{-}^{s+1})^{l+w+1-k}\|\chi_{B_l}\|_{L^{q'(\cdot)}(\mathbb{R}^{n})}\|\chi_{B_k}\|_{L^{q(\cdot)}(\mathbb{R}^{n})}\\
&\displaystyle\leq
Cb^{-l\alpha(0)-l+k}b^{(l-k)\delta_2}(b\lambda_{-}^{s+1})^{l+w+1-k}.

\end{array}\eqno(5.4)$$
Let $z=\lambda_{-}^{-(s+1)}b^{\alpha(0)-\delta_2}$. Then by $z<1$ we obtain
$$\begin{array}{rl}
\displaystyle V_1&\displaystyle\leq
C\sum_{k=-\infty}^{-1}b^{k\alpha(0)p(1+\varepsilon)}\left(\sum_{l=-\infty}^{k-w-1}|\lambda_l|b^{-l\alpha(0)-l+k}b^{(l-k)\delta_2}(b\lambda_{-}^{s+1})^{l+w+1-k}\right)^{p(1+\varepsilon)}\\
&\displaystyle\leq
C\sum_{k=-\infty}^{-1}\left(\sum_{l=-\infty}^{k-w-1}|\lambda_l|z^{k-l}\right)^{p(1+\varepsilon)}\\
&\displaystyle\leq C\sum_{k=-\infty}^{-1}z^{kp(1+\varepsilon)/2}\left(\sum_{l=-\infty}^{k-w-1}|\lambda_l|^{p(1+\varepsilon)}z^{-lp(1+\varepsilon)/2}\right)\left(\sum_{l=-\infty}^{k-w-1}z^{(k-l)
[p(1+\varepsilon)]'/2}\right)^{\frac{p(1+\varepsilon)}{[p(1+\varepsilon)]'}}\\
&\displaystyle\leq C\sum_{k=-\infty}^{-1}z^{kp(1+\varepsilon)/2}\left(\sum_{l=-\infty}^{k-w-1}|\lambda_l|^{p(1+\varepsilon)}z^{-lp(1+\varepsilon)/2}\right)\\
&\displaystyle\leq C\left(\sum_{l=-\infty}^{-w-2}|\lambda_l|^{p(1+\varepsilon)}z^{-lp(1+\varepsilon)/2}\right)\left(\sum_{k=l+w+1}^{-1}z^{kp(1+\varepsilon)/2}\right)\\
&\displaystyle\leq C\sum_{l=-\infty}^{-w-2}|\lambda_l|^{p(1+\varepsilon)}.

\end{array}$$
So we have
$$\sum_{k=-\infty}^{-1}\|b^{k\alpha(\cdot)}(M^0_Nf)\chi_k\|^{p(1+\varepsilon)}_{L^{q(\cdot)}(\mathbb{R}^{n})}\leq C\sum_{l=-\infty}^\infty|\lambda_l|^{p(1+\varepsilon)}.$$

When $k\geq 0$, we have
$$\begin{array}{rl}
\displaystyle
\sum_{k=0}^{\infty}\|b^{k\alpha(\cdot)}(M^0_Nf)\chi_k\|^{p(1+\varepsilon)}_{L^{q(\cdot)}(\mathbb{R}^{n})}
&\displaystyle\leq
\sum_{k=0}^{\infty}b^{k\alpha_\infty p(1+\varepsilon)}\left(\sum_{l=-\infty}^\infty|\lambda_l|\|(M^0_Na_l)\chi_k\|_{L^{q(\cdot)}(\mathbb{R}^{n})}\right)^{p(1+\varepsilon)}\\
&\displaystyle\leq C\sum_{k=0}^{\infty}b^{k\alpha_\infty p(1+\varepsilon)}\left(\sum_{l=-\infty}^{-1}|\lambda_l|\|(M^0_Na_l)\chi_k\|_{L^{q(\cdot)}(\mathbb{R}^{n})}\right)^{p(1+\varepsilon)}\\
&\displaystyle\hspace{3mm}+ C\sum_{k=0}^{\infty}b^{k\alpha_\infty p(1+\varepsilon)}\left(\sum_{l=0}^{k-w-1}|\lambda_l|\|(M^0_Na_l)\chi_k\|_{L^{q(\cdot)}(\mathbb{R}^{n})}\right)^{p(1+\varepsilon)}\\

\end{array}$$
$$\begin{array}{rl}
&\displaystyle\hspace{3mm}+C\sum_{k=0}^{\infty}b^{k\alpha_\infty p(1+\varepsilon)}\left(\sum_{l=k-w}^\infty|\lambda_l|\|(M^0_Na_l)\chi_k\|_{L^{q(\cdot)}(\mathbb{R}^{n})}\right)^{p(1+\varepsilon)}\\
&\displaystyle=W_1+W_2+W_3.

\end{array}$$
Using the $L^{q(\cdot)}(\mathbb{R}^{n})$ boundedness of the radial grand
maximal operator $M^0_N$ and the H\"{o}lder inequality, we have
$$\begin{array}{rl}
\displaystyle W_3&\displaystyle\leq
C\sum_{k=0}^{\infty}b^{k\alpha_\infty p(1+\varepsilon)}\left(\sum_{l=k-w}^\infty|\lambda_l|\|a_l\|_{L^{q(\cdot)}(\mathbb{R}^{n})}\right)^{p(1+\varepsilon)}\\
&\displaystyle\leq C\sum_{k=0}^{\infty}\left(\sum_{l=k-w}^\infty|\lambda_l|^{p(1+\varepsilon)}b^{(k-l)\alpha_\infty
p(1+\varepsilon)/2}\right)\left(\sum_{l=k-w}^\infty b^{(k-l)\alpha_\infty
[p(1+\varepsilon)]'/2}\right)^{\frac{p(1+\varepsilon)}{[p(1+\varepsilon)]'}}\\
&\displaystyle\leq C\sum_{l=-w}^\infty|\lambda_l|^{p(1+\varepsilon)}\left(\sum_{k=0}^{l+w}b^{(k-l)\alpha_\infty
p(1+\varepsilon)/2}\right)\\
&\displaystyle\leq C\sum_{l=-w}^\infty|\lambda_l|^{p(1+\varepsilon)}.

\end{array}$$
For $W_1$, let $z=\lambda_{-}^{-(s+1)}b^{\alpha^+-\delta_2}$. Similar to the estimate of $V_1$, by (4.1), $z<1$ and $\delta_2\leq\alpha(0),\,\alpha_\infty\leq \alpha^+$, we have
$$\begin{array}{rl}
\displaystyle W_1&\displaystyle\leq
C\sum_{k=0}^{\infty}b^{k\alpha_\infty p(1+\varepsilon)}\left(\sum_{l=-\infty}^{-1}|\lambda_l|b^{-l\alpha(0)-l+k}b^{(l-k)\delta_2}(b\lambda_{-}^{s+1})^{l+w+1-k}\right)^{p(1+\varepsilon)}\\
&\displaystyle\leq
C\sum_{k=0}^{\infty}\left(\sum_{l=-\infty}^{-1}|\lambda_l|z^{k-l}\right)^{p(1+\varepsilon)}\\
&\displaystyle\leq C\sum_{k=0}^{\infty}z^{kp(1+\varepsilon)}\left(\sum_{l=-\infty}^{-1}|\lambda_l|^{p(1+\varepsilon)}z^{-lp(1+\varepsilon)/2}\right)\left(\sum_{l=-\infty}^{-1}z^{-l
[p(1+\varepsilon)]'/2}\right)^{\frac{p(1+\varepsilon)}{[p(1+\varepsilon)]'}}\\
&\displaystyle\leq C\sum_{l=-\infty}^{-1}|\lambda_l|^{p(1+\varepsilon)}.

\end{array}$$
For $W_2$, let $z=\lambda_{-}^{-(s+1)}b^{\alpha_\infty-\delta_2}$. Similar to the estimate of $V_1$, by (4.1) and $z<1$ we have
$$\begin{array}{rl}
\displaystyle W_2&\displaystyle\leq
C\sum_{k=0}^{\infty}b^{k\alpha_\infty p(1+\varepsilon)}\left(\sum_{l=0}^{k-w-1}|\lambda_l|b^{-l\alpha_\infty-l+k}b^{(l-k)\delta_2}(b\lambda_{-}^{s+1})^{l+w+1-k}\right)^{p(1+\varepsilon)}\\
&\displaystyle\leq
C\sum_{k=0}^{\infty}\left(\sum_{l=0}^{k-w-1}|\lambda_l|z^{k-l}\right)^{p(1+\varepsilon)}\\
&\displaystyle\leq C\sum_{k=0}^{\infty}z^{kp(1+\varepsilon)/2}\left(\sum_{l=0}^{k-w-1}|\lambda_l|^{p(1+\varepsilon)}z^{-lp(1+\varepsilon)/2}\right)\left(\sum_{l=0}^{k-w-1}z^{(k-l)
[p(1+\varepsilon)]'/2}\right)^{\frac{p(1+\varepsilon)}{[p(1+\varepsilon)]'}}\\
&\displaystyle\leq C\sum_{k=0}^{\infty}z^{kp(1+\varepsilon)/2}\left(\sum_{l=0}^{k-w-1}|\lambda_l|^{p(1+\varepsilon)}z^{-lp(1+\varepsilon)/2}\right)\\
&\displaystyle\leq C\left(\sum_{l=0}^{\infty}|\lambda_l|^{p(1+\varepsilon)}z^{-lp(1+\varepsilon)/2}\right)\left(\sum_{k=l+w+1}^{\infty}z^{kp(1+\varepsilon)/2}\right)\\

\end{array}$$
$$\begin{array}{rl}
&\displaystyle\leq C\sum_{l=0}^{\infty}|\lambda_l|^{p(1+\varepsilon)}.

\end{array}$$
So we have
$$\sum_{k=0}^{\infty}\|b^{k\alpha(\cdot)}(M^0_Nf)\chi_k\|^{p(1+\varepsilon)}_{L^{q(\cdot)}(\mathbb{R}^{n})}\leq C\sum_{l=-\infty}^\infty|\lambda_l|^{p(1+\varepsilon)}.$$

Thus we can get
$$\begin{array}{rl}
\displaystyle\|f\|^{p(1+\varepsilon)}_{H\dot{K}^{\alpha(\cdot),p)}_{q(\cdot),\theta}(A;\mathbb{R}^{n})}&\displaystyle\leq C
\|M^0_Nf\|^{p(1+\varepsilon)}_{\dot{K}^{\alpha(\cdot),p),\theta}_{q(\cdot)}(A;\mathbb{R}^{n})}\\
&\displaystyle=C\sup_{\varepsilon>0}\varepsilon^\theta\sum_{k=-\infty}^\infty\|b^{k\alpha(\cdot)}(M^0_Nf)\chi_k\|^{p(1+\varepsilon)}_{L^{q(\cdot)}(\mathbb{R}^{n})}\\
&\displaystyle\leq
C\sup_{\varepsilon>0}\varepsilon^\theta\sum_{l=-\infty}^\infty|\lambda_l|^{p(1+\varepsilon)}<\infty.

\end{array}$$

This finishes the proof of Theorem 5.1.

\noindent {\bf Remark 5.2}\quad For the sufficiency of Theorem 5.1, if we
change the condition $\delta_2\leq\alpha(0), \alpha_\infty<\delta_2+\ln\lambda_{-}/\ln b$ to $\delta_2\leq\alpha(0), \alpha_\infty<\infty$, then the conclusion of Theorem 5.1 is also true.
\medskip

As an application of the atomic decomposition theorems, we shall
extend Theorem 3.3 to the case of $\alpha\geq n\delta_2$. And the following condition (5.5) is necessary for our discussion on the linear operator $T$ on anisotropic grand Herz-type Hardy spaces with variable exponents:
$$f=\sum_{i\in\mathbb{N}}\lambda_ia_i\quad\mathrm{in}\quad \mathcal{S}'\quad\Rightarrow\quad Tf=\sum_{i\in\mathbb{N}}\lambda_iTa_i \quad\mathrm{in}\quad \mathcal{S}'.\eqno(5.5)$$

\noindent{\bf Theorem 5.2}\quad Let $\alpha(\cdot)\in L^\infty(\mathbb{R}^{n})\cap\mathcal{P}_0(\mathbb{R}^{n})\cap\mathcal{P}_\infty(\mathbb{R}^{n})$, $\delta_2\leq\alpha(0), \alpha_\infty<\delta_2+\ln\lambda_{-}/\ln b$,
$1\leq p<\infty$ and $q(\cdot)\in \mathcal{B}(\mathbb{R}^{n})$. If a linear operator $T$ satisfies (5.5) for every central atomic decomposition, is bounded on $L^{q(\cdot)}(\mathbb{R}^{n})$, and for any $f\in L^{q(\cdot)}(\mathbb{R}^{n})$ with compact support $B_k$, $\int_{\mathbb{R}^{n}}f(x)dx=0$, $T$ satisfies the following size condition:
$$|Tf(x)|\leq C\frac{\|f\|_{L^1}}{(\rho(x))^2}, \quad \mathrm{if} \quad \inf_{y\in{\rm
supp}f}\rho(x-y)\geq b^{-w}\left(1-\frac{1}{b}\right)\rho(x),\eqno(5.6)$$
then $T$ can be extended to be a bounded operator from
$H\dot{K}^{\alpha(\cdot),p)}_{q(\cdot),\theta}(A;\mathbb{R}^{n})$ to
$\dot{K}^{\alpha(\cdot),p)}_{q(\cdot),\theta}(A;\mathbb{R}^{n})$ (or bounded from
$HK^{\alpha(\cdot),p)}_{q(\cdot),\theta}(A;\mathbb{R}^{n})$ to
$K^{\alpha(\cdot),p)}_{q(\cdot),\theta}(A;\mathbb{R}^{n})$).

\noindent{\bf Proof}\quad It suffices to prove homogeneous case.
Suppose $f\in H\dot{K}^{\alpha(\cdot),p)}_{q(\cdot),\theta}(A;\mathbb{R}^{n})$. By
Theorem 5.1, $\displaystyle f=\sum_{j=-\infty}^\infty\lambda_ja_j$
in the sense of $\mathcal {S'}(\mathbb{R}^{n})$, where each $a_j$ is
a central $(\alpha(\cdot), q(\cdot))$-atom with support contained in $B_j$
and
$$\|f\|_{H\dot{K}^{\alpha(\cdot),p)}_{q(\cdot),\theta}(A;\mathbb{R}^{n})}\approx\inf\displaystyle\sup_{\varepsilon>0}\left(\varepsilon^\theta\sum_{j=-\infty}^\infty|\lambda_j|^{p(1+\varepsilon)}\right)^{\frac{1}{p(1+\varepsilon)}},$$
where the infimum is taken over all above decompositions of $f$. Therefore, we get
$$\begin{array}{rl}
\displaystyle \|Tf\|^{p(1+\varepsilon)}_{\dot{K}^{\alpha(\cdot),p)}_{q(\cdot),\theta}(A;\mathbb{R}^{n})}&\displaystyle=\sup_{\varepsilon>0}\varepsilon^\theta\sum_{k=-\infty}^\infty
\|b^{k\alpha(\cdot)}(Tf)\chi_k\|^{p(1+\varepsilon)}_{L^{q(\cdot)}(\mathbb{R}^{n})}\\
&\displaystyle\leq\sup_{\varepsilon>0}\varepsilon^\theta\sum_{k=-\infty}^{-1}b^{k\alpha(0)p(1+\varepsilon)}\|(Tf)\chi_k\|^{p(1+\varepsilon)}_{L^{q(\cdot)}(\mathbb{R}^{n})}\\
&\displaystyle\hspace{3mm}+\sup_{\varepsilon>0}\varepsilon^\theta\sum_{k=0}^\infty b^{k\alpha_\infty p(1+\varepsilon)}\|(Tf)\chi_k\|^{p(1+\varepsilon)}_{L^{q(\cdot)}(\mathbb{R}^{n})}\\
&\displaystyle\leq C\bigg[\sup_{\varepsilon>0}\varepsilon^\theta\sum_{k=-\infty}^{-1}b^{k\alpha(0)p(1+\varepsilon)}\bigg(\sum_{j=-\infty}^{k-w-1}|\lambda_j|\|(Ta_j)\chi_k\|_{L^{q(\cdot)}
(\mathbb{R}^{n})}\bigg)^{p(1+\varepsilon)}\\
&\displaystyle\hspace{3mm}+\sup_{\varepsilon>0}\varepsilon^\theta\sum_{k=-\infty}^{-1}b^{k\alpha(0)p(1+\varepsilon)}\bigg(\sum_{j=k-w}^{-1}|\lambda_j|\|(Ta_j)\chi_k\|_{L^{q(\cdot)}
(\mathbb{R}^{n})}\bigg)^{p(1+\varepsilon)}\\
&\displaystyle\hspace{3mm}+\sup_{\varepsilon>0}\varepsilon^\theta\sum_{k=-\infty}^{-1}b^{k\alpha(0)p(1+\varepsilon)}\bigg(\sum_{j=0}^{\infty}|\lambda_j|\|(Ta_j)\chi_k\|_{L^{q(\cdot)}
(\mathbb{R}^{n})}\bigg)^{p(1+\varepsilon)}\\
&\displaystyle\hspace{3mm}+\sup_{\varepsilon>0}\varepsilon^\theta\sum_{k=0}^\infty b^{k\alpha_\infty p(1+\varepsilon)}\bigg(\sum_{j=-\infty}^{-1}|\lambda_j|\|(Ta_j)\chi_k\|_{L^{q(\cdot)}
(\mathbb{R}^{n})}\bigg)^{p(1+\varepsilon)}\\
&\displaystyle\hspace{3mm}+\sup_{\varepsilon>0}\varepsilon^\theta\sum_{k=0}^\infty b^{k\alpha_\infty p(1+\varepsilon)}\bigg(\sum_{j=0}^{k-w-1}|\lambda_j|\|(Ta_j)\chi_k\|_{L^{q(\cdot)}
(\mathbb{R}^{n})}\bigg)^{p(1+\varepsilon)}\\
&\displaystyle\hspace{3mm}+\sup_{\varepsilon>0}\varepsilon^\theta\sum_{k=0}^\infty b^{k\alpha_\infty p(1+\varepsilon)}\bigg(\sum_{j=k-w}^{\infty}|\lambda_j|\|(Ta_j)\chi_k\|_{L^{q(\cdot)}
(\mathbb{R}^{n})}\bigg)^{p(1+\varepsilon)}\bigg]\\
&\displaystyle= C(U_1+U_2+U_3+U_4+U_5+U_6).
\end{array}\eqno(5.7)$$

Let us first estimate $U_1$. By (5.6) and the generalized H\"{o}lder
inequality, we get
$$|Ta_j(x)|\leq C\frac{b^j\|a_j\|_{L^1}}{(\rho(x))^2}\leq Cb^{j+2-2k}\|a_j\|_{L^{q(\cdot)}(\mathbb{R}^{n})}\|\chi_{B_j}\|_{L^{q'(\cdot)}(\mathbb{R}^{n})}.$$
So by Lemma 2.2 and Lemma 2.3, we have
$$\begin{array}{rl}
\displaystyle \|(Ta_j)\chi_k\|_{L^{q(\cdot)}(\mathbb{R}^{n})}&\displaystyle\lesssim b^{j+2-2k}\|a_j\|_{L^{q(\cdot)}(\mathbb{R}^{n})}\|\chi_{B_j}\|_{L^{q'(\cdot)}(\mathbb{R}^{n})}\|\chi_{B_k}\|_{L^{q(\cdot)}(\mathbb{R}^{n})}\\
&\displaystyle\lesssim b^{j+2-2k}\|a_j\|_{L^{q(\cdot)}(\mathbb{R}^{n})}\big(|B_k|\|\chi_{B_k}\|^{-1}_{L^{q'(\cdot)}(\mathbb{R}^{n})}\big)\|\chi_{B_j}\|_{L^{q'(\cdot)}(\mathbb{R}^{n})}\\
&\displaystyle\lesssim b^{j+2-k}b^{\delta_2(j-k)}\|a_j\|_{L^{q(\cdot)}(\mathbb{R}^{n})}.

\end{array}\eqno(5.8)$$
Since $b>1, \alpha(0)<\delta_2+1$, by
(5.8) and the H\"{o}lder inequality, we have
$$\begin{array}{rl}
\displaystyle U_1&\displaystyle\lesssim\sup_{\varepsilon>0}\varepsilon^\theta\sum_{k=-\infty}^{-1}b^{k\alpha(0) p(1+\varepsilon)}\bigg(\sum_{j=-\infty}^{k-w-1}|\lambda_j|b^{j+2-k}b^{\delta_2(j-k)-j\alpha(0)}\bigg)^{p(1+\varepsilon)}\\
&\displaystyle\lesssim\sup_{\varepsilon>0}\varepsilon^\theta\sum_{k=-\infty}^{-1}\bigg(\sum_{j=-\infty}^{k-w-1}|\lambda_j|b^{(j-k)(\delta_2-\alpha(0)+1)}\bigg)^{p(1+\varepsilon)}\\
&\displaystyle\lesssim\sup_{\varepsilon>0}\varepsilon^\theta\sum_{k=-\infty}^{-1}\bigg(\sum_{j=-\infty}^{k-w-1}|\lambda_j|^{p(1+\varepsilon)}b^{(j-k)[\delta_2-\alpha(0)+1]p(1+\varepsilon)/2}\bigg)
\bigg(\sum_{j=-\infty}^{k-w-1}b^{(j-k)[\delta_2-\alpha(0)+1][p(1+\varepsilon)]'/2}\bigg)^{\frac{p(1+\varepsilon)}{[p(1+\varepsilon)]'}}\\
&\displaystyle\lesssim\sup_{\varepsilon>0}\varepsilon^\theta\sum_{k=-\infty}^{-1}\bigg(\sum_{j=-\infty}^{k-w-1}|\lambda_j|^{p(1+\varepsilon)}b^{(j-k)[\delta_2-\alpha(0)+1]p(1+\varepsilon)/2}\bigg)\\

\end{array}$$
$$\begin{array}{rl}
&\displaystyle\lesssim\sup_{\varepsilon>0}\varepsilon^\theta\sum_{j=-\infty}^{-w-2}|\lambda_j|^{p(1+\varepsilon)}\sum_{k=j+w+1}^{-1}b^{(j-k)[\delta_2-\alpha(0)+1]p(1+\varepsilon)/2}\\
&\displaystyle\lesssim\sup_{\varepsilon>0}\varepsilon^\theta\sum_{j=-\infty}^{-w-2}|\lambda_j|^{p(1+\varepsilon)}.

\end{array}\eqno(5.9)$$

Let us now estimate $U_2$. By
$L^{q(\cdot)}(\mathbb{R}^{n})$ boundedness of $T$ and the H\"{o}lder
inequality, we have
$$\begin{array}{rl}
\displaystyle U_2&\displaystyle\lesssim\sup_{\varepsilon>0}\varepsilon^\theta\sum_{k=-\infty}^{-1}b^{k\alpha(0) p(1+\varepsilon)}\bigg(\sum_{j=k-w}^{-1}|\lambda_j|\|a_j\|_{L^{q(\cdot)}(\mathbb{R}^{n})}\bigg)^{p(1+\varepsilon)}\\
&\displaystyle\lesssim\sup_{\varepsilon>0}\varepsilon^\theta\sum_{k=-\infty}^{-1}\bigg(\sum_{j=k-w}^{-1}|\lambda_j|b^{(k-j)\alpha(0)}\bigg)^{p(1+\varepsilon)}\\
&\displaystyle\lesssim\sup_{\varepsilon>0}\varepsilon^\theta\sum_{k=-\infty}^{-1}\bigg(\sum_{j=k-w}^{-1}|\lambda_j|^{p(1+\varepsilon)}b^{(k-j)\alpha(0)p(1+\varepsilon)/2}\bigg)\bigg(\sum_{j=k-w}^{-1}b^{(k-j)\alpha(0)[p(1+\varepsilon)]'/2}\bigg)^{p(1+\varepsilon)/[p(1+\varepsilon)]'}\\
&\displaystyle\lesssim\sup_{\varepsilon>0}\varepsilon^\theta\sum_{j=-\infty}^{-1}|\lambda_j|^{p(1+\varepsilon)}\sum_{k=-\infty}^{j+w}b^{(k-j)\alpha(0)p(1+\varepsilon)/2}\\
&\displaystyle\lesssim\sup_{\varepsilon>0}\varepsilon^\theta\sum_{j=-\infty}^{-1}|\lambda_j|^{p(1+\varepsilon)}.

\end{array}\eqno(5.10)$$

Let us now estimate $U_3$. By
$L^{q(\cdot)}(\mathbb{R}^{n})$ boundedness of $T$ and the H\"{o}lder
inequality, we have
$$\begin{array}{rl}
\displaystyle U_3&\displaystyle\lesssim\sup_{\varepsilon>0}\varepsilon^\theta\sum_{k=-\infty}^{-1}b^{k\alpha(0) p(1+\varepsilon)}\bigg(\sum_{j=0}^{\infty}|\lambda_j|\|a_j\|_{L^{q(\cdot)}(\mathbb{R}^{n})}\bigg)^{p(1+\varepsilon)}\\
&\displaystyle\lesssim\sup_{\varepsilon>0}\varepsilon^\theta\sum_{k=-\infty}^{-1}b^{k\alpha(0) p(1+\varepsilon)}\bigg(\sum_{j=0}^{\infty}|\lambda_j|b^{-j\alpha_\infty}\bigg)^{p(1+\varepsilon)}\\
&\displaystyle\lesssim\sup_{\varepsilon>0}\varepsilon^\theta\bigg(\sum_{j=0}^{\infty}|\lambda_j|^{p(1+\varepsilon)}b^{-j\alpha_\infty p(1+\varepsilon)/2}\bigg)\bigg(\sum_{j=0}^{\infty}b^{-j\alpha_\infty [p(1+\varepsilon)]'/2}\bigg)^{p(1+\varepsilon)/[p(1+\varepsilon)]'}\\
&\displaystyle\lesssim\sup_{\varepsilon>0}\varepsilon^\theta\sum_{j=0}^{\infty}|\lambda_j|^{p(1+\varepsilon)}.

\end{array}\eqno(5.11)$$

For $U_4$. By
(5.8), $\delta_2\leq\alpha(0),\alpha_\infty<\delta_2+1$ and the H\"{o}lder inequality, we have
$$\begin{array}{rl}
\displaystyle U_4&\displaystyle\lesssim\sup_{\varepsilon>0}\varepsilon^\theta\sum_{k=0}^\infty b^{k\alpha_\infty p(1+\varepsilon)}\bigg(\sum_{j=-\infty}^{-1}|\lambda_j|b^{j+2-k}b^{\delta_2(j-k)-j\alpha(0)}\bigg)^{p(1+\varepsilon)}\\
&\displaystyle\lesssim\sup_{\varepsilon>0}\varepsilon^\theta\sum_{k=0}^\infty b^{k(\alpha_\infty-\delta_2-1) p(1+\varepsilon)}\bigg(\sum_{j=-\infty}^{-1}|\lambda_j|b^{j[1+\delta_2-\alpha(0)]}\bigg)^{p(1+\varepsilon)}\\
&\displaystyle\lesssim\sup_{\varepsilon>0}\varepsilon^\theta\bigg(\sum_{j=-\infty}^{-1}|\lambda_j|^{p(1+\varepsilon)}b^{j[1+\delta_2-\alpha(0)]p(1+\varepsilon)/2}\bigg)\bigg(\sum_{j=-\infty}^{-1}b^{j[1+\delta_2-\alpha(0)][p(1+\varepsilon)]'/2}\bigg)^{p(1+\varepsilon)/[p(1+\varepsilon)]'}\\
&\displaystyle\lesssim\sup_{\varepsilon>0}\varepsilon^\theta\sum_{j=-\infty}^{-1}|\lambda_j|^{p(1+\varepsilon)}b^{j[1+\delta_2-\alpha(0)]p(1+\varepsilon)/2}\\

\end{array}$$
$$\begin{array}{rl}
&\displaystyle\lesssim\sup_{\varepsilon>0}\varepsilon^\theta\sum_{j=-\infty}^{-1}|\lambda_j|^{p(1+\varepsilon)}.

\end{array}\eqno(5.12)$$

For $U_5$. By
(5.8), $\delta_2\leq\alpha(0),\alpha_\infty<\delta_2+1$ and the H\"{o}lder inequality, we have
$$\begin{array}{rl}
\displaystyle U_5&\displaystyle\lesssim\sup_{\varepsilon>0}\varepsilon^\theta\sum_{k=0}^\infty b^{k\alpha_\infty p(1+\varepsilon)}\bigg(\sum_{j=0}^{k-w-1}|\lambda_j|b^{j+2-k}b^{\delta_2(j-k)-j\alpha_\infty}\bigg)^{p(1+\varepsilon)}\\
&\displaystyle\lesssim\sup_{\varepsilon>0}\varepsilon^\theta\sum_{k=0}^\infty \bigg(\sum_{j=0}^{k-w-1}|\lambda_j|b^{(j-k)(1+\delta_2-\alpha_\infty)}\bigg)^{p(1+\varepsilon)}\\
&\displaystyle\lesssim\sup_{\varepsilon>0}\varepsilon^\theta\sum_{k=0}^\infty\bigg(\sum_{j=0}^{k-w-1}|\lambda_j|^{p(1+\varepsilon)}b^{(j-k)(1+\delta_2-\alpha_\infty)p(1+\varepsilon)/2}\bigg)\bigg(\sum_{j=0}^{k-w-1}b^{(j-k)(1+\delta_2-\alpha_\infty)[p(1+\varepsilon)]'/2}\bigg)^{\frac{p(1+\varepsilon)}{[p(1+\varepsilon)]'}}\\
&\displaystyle\lesssim\sup_{\varepsilon>0}\varepsilon^\theta\sum_{k=0}^\infty\sum_{j=0}^{k-w-1}|\lambda_j|^{p(1+\varepsilon)}b^{(j-k)(1+\delta_2-\alpha_\infty)p(1+\varepsilon)/2}\\
&\displaystyle\lesssim\sup_{\varepsilon>0}\varepsilon^\theta\sum_{j=0}^\infty|\lambda_j|^{p(1+\varepsilon)}\sum_{k=j+w+1}^\infty b^{(j-k)(1+\delta_2-\alpha_\infty)p(1+\varepsilon)/2}\\
&\displaystyle\lesssim\sup_{\varepsilon>0}\varepsilon^\theta\sum_{j=0}^\infty|\lambda_j|^{p(1+\varepsilon)}.

\end{array}\eqno(5.13)$$

Let us now estimate $U_6$. By
$L^{q(\cdot)}(\mathbb{R}^{n})$ boundedness of $T$ and the H\"{o}lder
inequality, we have
$$\begin{array}{rl}
\displaystyle U_6&\displaystyle=\sup_{\varepsilon>0}\varepsilon^\theta\sum_{k=0}^\infty b^{k\alpha_\infty p(1+\varepsilon)}\bigg(\sum_{j=k-w}^{\infty}|\lambda_j|\|(Ta_j)\chi_k\|_{L^{q(\cdot)}(\mathbb{R}^{n})}\bigg)^{p(1+\varepsilon)}\\
&\displaystyle\lesssim\sup_{\varepsilon>0}\varepsilon^\theta\sum_{k=0}^\infty b^{k\alpha_\infty p(1+\varepsilon)}\bigg(\sum_{j=k-w}^{\infty}|\lambda_j|\|a_j\|_{L^{q(\cdot)}(\mathbb{R}^{n})}\bigg)^{p(1+\varepsilon)}\\
&\displaystyle\lesssim\sup_{\varepsilon>0}\varepsilon^\theta\sum_{k=0}^\infty b^{k\alpha_\infty p(1+\varepsilon)}\bigg(\sum_{j=k-w}^\infty|\lambda_j|b^{-j\alpha_\infty}\bigg)^{p(1+\varepsilon)}\\
&\displaystyle\lesssim\sup_{\varepsilon>0}\varepsilon^\theta\sum_{k=0}^\infty b^{k\alpha_\infty p(1+\varepsilon)}\bigg(\sum_{j=k-w}^\infty|\lambda_j|^{p(1+\varepsilon)}b^{-j\alpha_\infty p(1+\varepsilon)/2}\bigg)\bigg(\sum_{j=k-w}^\infty b^{-j\alpha_\infty [p(1+\varepsilon)]'/2}\bigg)^{p(1+\varepsilon)/[p(1+\varepsilon)]'}\\
&\displaystyle\lesssim\sup_{\varepsilon>0}\varepsilon^\theta\sum_{k=0}^\infty b^{k\alpha_\infty p(1+\varepsilon)/2}\bigg(\sum_{j=k-w}^\infty|\lambda_j|^{p(1+\varepsilon)}b^{-j\alpha_\infty p(1+\varepsilon)/2}\bigg)\\
&\displaystyle\lesssim\sup_{\varepsilon>0}\varepsilon^\theta\sum_{j=-w}^\infty|\lambda_j|^{p(1+\varepsilon)}\sum_{k=0}^{j+w}b^{(k-j)\alpha_\infty p(1+\varepsilon)/2}\\
&\displaystyle\lesssim\sup_{\varepsilon>0}\varepsilon^\theta\sum_{j=-w}^\infty|\lambda_j|^{p(1+\varepsilon)}.

\end{array}\eqno(5.14)$$
Combining (5.7), (5.9)-(5.14), we have
$$\|Tf\|_{\dot{K}^{\alpha(\cdot),p)}_{q(\cdot),\theta}(A;\mathbb{R}^{n})}\leq C\|f\|_{H\dot{K}^{\alpha(\cdot),p)}_{q(\cdot),\theta}(A;\mathbb{R}^{n})}.$$

Thus, the proof of Theorem 5.2 is completed.

\vskip5mm
\noindent {\bf Acknowledgement} The authors are very grateful to the referees
for their valuable comments. This work was supported by Shandong Provincial Natural Science Foundation (Grant No. ZR2022MA054) and National Natural Science Foundation of China (Grant Nos. 12271267, 11926343 and 11926342).

\vskip5mm
\centerline{\Large\bf References} \vspace{0.4cm} \def%
\hang{\hangindent\parindent} \def\textindent#1{\indent\llap{#1\enspace}%
\ignorespaces} \def\re{\par\hang\textindent}

\re{[1]} A. Almeida and D. Drihem, Maximal, potential and singular type operators on Herz spaces with variable exponents, {\it J. Math. Anal. Appl.,} 394(2012), 781-795.

\re{[2]} M. Bownik, Anisotropic Hardy spaces and wavelets, \textit{Mem. Amer. Math. Soc.} \textbf{164}
(2003), 122 pp.

\re{[3]} A. P. Calder\'{o}n and A. Torchinsky, Parabolic maximal functions associated with a distribution, \textit{Adv. Math.} \textbf{16}
(1975), 1-64.

\re{[4]} A. P. Calder\'{o}n and A. Torchinsky, Parabolic maximal functions associated with a distribution II, \textit{Adv. Math.} \textbf{24}
(1977), 101-171.

\re{[5]} C. Capone, M. Formica and R. Giova, R, Grand Lebesgue spaces with respect to measurable functions, {\it Nonlinear Anal.,} {\bf 85}(2013), 125-131.

\re{[6]} Y. Chen, S. Levin and M. Rao, Variable exponent, linear growth functionals in image restoration, {\it SIAM J. Appl. Math.,} {\bf 66}(2006), 1383-1406.

\re{[7]} D. Cruz-Uribe and A. Fiorenza, Variable Lebesgue Spaces: Foundations and Harmonic Analysis(Applied and Numerical
Harmonic Analysis), Springer, Heidelberg, 2013.

\re{[8]} L. Diening, P. Harjulehto, P. H\"{a}st\"{o} and M.
R\r{u}\v{z}i\v{c}ka, {\it Lebesgue and Sobolev spaces with variable
exponents,} Lecture Notes in Math., vol. {\bf 2017},
Springer-Verlag, Berlin, 2011.

\re{[9]} Y. Ding, S. Lan and S. Lu, New Hardy spaces associated with some anisotropic
Herz spaces and their applications, {\it Acta Math. Sin.,} {\bf 24}(2008), 1449-1470.

\re{[10]} D. Drihem and F. Seghiri, Notes on the Herz-type Hardy spaces of variable smoothness and integrability, {\it Math. Inequal. Appl.,} {\bf 19}(2016), 145-165.

\re{[11]} C. Fefferman and E. M. Stein, $H^p$ spaces of several variables, \textit{Acta Math.} \textbf{129}
(1972), 137-193.

\re{[12]} A. Fiorenza, Duality and reflexivity in grand Lebesgue spaces, {\it Collect. Math.,} {\bf 51}(2000), 131-148.

\re{[13]} A. Fiorenza, A. Mercaldo and J. Rakotoson, Regularity and comparison results in grand Sobolev spaces for parabolic equations with measure data, {\it Appl. Math. Lett.,} {\bf 14}(2001), 979-981.

\re{[14]} Z. Fu, Z. Liu, S. Lu and H. Wang, Characterization for commutators of $n$-dimensional fractional Hardy operators, {\it  Sci. China Ser. A-Math.,} {\bf 50}(2007), 1418-1426.

\re{[15]} P. Harjulehto, P. H\"{a}st\"{o}, \'{U}. V. L\^{e} and M.
Nuortio, Overview of differential equations with non-standard
growth, {\it Nonlinear Anal.,} {\bf 72} (2010), 4551-4574.

\re{[16]} C. Herz, Lipshitz spaces and Bernstein's theorem on absolutely convergent Fourier transforms, {\it J. Math. Mech.,} {\bf 18}(1968), 283-324.

\re{[17]} T. Iwaniec and C. Sbordone, On the integrability of the Jacobian under minimal
hypotheses, {\it Arch. Rational Mech. Anal.,} {\bf 119}(1992), 129-143.

\re{[18]} M. Izuki, Boundedness of sublinear operators on Herz spaces
with variable exponent and application to wavelet characterization,
{\it Anal. Math.,} {\bf 36} (2010), 33-50.

\re{[19]} M. Izuki, Fractional integrals on Herz-Morrey spaces
with variable exponent, {\it Hiroshima Math. J.,} {\bf 40}(2010), 343-355.

\re{[20]} V. Kokilashvili, A. Meskhi and H. Rafeiro, Boundedness of commutators of singular and potential operators in generalized grand Morrey spaces and some applications, {\it Studia Math.,} {\bf 217}(2013), 159-178.

\re{[21]} Y. Komori, Weak type estimates for Calder\'{o}n-Zygmund operators on Herz spaces at critical indexes,
{\it Math. Nachr.,} {\bf 259}(2003), 42-50.

\re{[22]} O. Kov\'{a}\v{c}ik and J. R\'{a}kosn\'{i}k, On spaces
$L^{p(x)}$ and $W^{k,p(x)}$, {\it Czechoslovak Math. J.,} {\bf
41}(1991), 592-618.

\re{[23]} S. Lan, The decompositon of anisotropic Herz spaces and its applications, \textit{Panamer. Math. J.} \textbf{16}
(2006), 61-75.

\re{[24]} Z. Liu, Boundedness of commutators of fractional integration on Herz-type spaces. {\it Acta Math. Sci. Ser. B Engl. Ed.,} {\bf 20}(2000), 461-470.

\re{[25]} J. Liu, D. Yang and W. Yuan, Anisotropic Hardy-Lorentz spaces and their
applications, {\it Sci. China Math.,} {\bf 59}(2016), 1669-1720.

\re{[26]} J. Liu, D. Yang and W. Yuan, Anisotropic variable Hardy-Lorentz spaces and their real
interpolation, {\it J. Math. Anal. Appl.,} {\bf 456}(2017), 356-393.

\re{[27]} S. Lu and L. Xu, Boundedness of rough singular integral operators on the
homogeneous Morrey-Herz spaces, {\it Hokkaido Math. J.,} {\bf 34} (2005), 299-314.

\re{[28]} S. Lu and D. Yang, The local versions of
$H^p(\mathbb{R}^{n})$ spaces at the origin, {\it Studia Math.,} {\bf
116} (1995), 103-131.

\re{[29]} S. Lu and D. Yang, The weighted Herz-type Hardy space and
its applications, {\it Sci. China Ser. A,} {\bf 38} (1995), 662-673.

\re{[30]} H. Nafis, H. Rafeiro and M. Zaighum, A note on the boundedness of sublinear operators on grand variable Herz spaces, {\it J. Inequal. Appl.,} {\bf 2020}(2020), 13 pp.

\re{[31]} W. Orlicz, \"{U}ber konjugierte Exponentenfolgen, {\it Studia Math.,} {\bf 3}(1931), 200-211.

\re{[32]} H. Rafeiro, S. Samko and S. Umarkhadzhiev, Grand Lebesgue sequence spaces, {\it Georgian Math. J.,} {\bf 25}(2018), 291-302.

\re{[33]} J. Ruan, D. Fan and Q Wu, Weighted Herz space estimates for Hausdorff operators on the Heisenberg group, {\it Banach J. Math. Anal.,} {\bf 11}(2017), 513-535.

\re{[34]} M. R\r{u}\v{z}i\v{c}ka, Electrorheological fluids: modeling and mathematical theory, Springer, Berlin, 2000.

\re{[35]} F. Shabbir and M. Zaighum, On the boundedness of sublinear operators on grand Herz-Hardy spaces
with variable exponent, {\it Mediterr. J. Math.,} {\bf 21}(2024), 20 pp.

\re{[36]} H. Wang, Anisotropic Herz spaces with variable exponents, {\it Commun. Math. Anal.,} {\bf 18}(2015), 1-14.

\re{[37]} H. Wang and Z. Liu, The Herz-type Hardy spaces with variable
exponent and their applications, {\it Taiwanese J. Math.,} {\bf
16}(2012), 1363-1389.

\re{[38]} H. Wang and Y. Wu, Anisotropic Herz-Morrey spaces with variable exponents, {\it Khayyam J. Math.,} {\bf 2}(2016), 177-187.

\re{[39]} H. Zhao and J. Zhou, Anisotropic Herz-type Hardy spaces with variable
exponent and their applications, {\it Acta Math. Hungar.,} {\bf 156}(2018), 309-335.

\bigskip

\medskip

\noindent\author{Hongbin \uppercase{Wang}}\\
    {School of Mathematics and Statistics, Shandong University of Technology, Zibo, Shandong, 255049, China\\
    E-mail\,$:$ hbwang\_2006@163.com}
\medskip

\noindent\author{Zongguang \uppercase{Liu}}\\
    {Department of Mathematics, China University of Mining and Technology(Beijing), Beijing, 100083, China\\
    E-mail\,$:$ liuzg@cumtb.edu.cn}

\end{document}